\documentclass[11pt]{article}
\usepackage{amsfonts}
\usepackage{color}
\newtheorem{theo}{Theorem}[section]
\newtheorem{lem}[theo]{Lemma}
\newtheorem{cor}[theo]{Corollary}

\newcommand{\mysection}[1]{\section{#1} \setcounter{equation}{0}}
\newcommand{\proof}{{\sc Proof.} \quad}
\newcommand{\proofc}{{\sc Proof} \ }
\newcommand{\be}{\begin{equation} \label}
\newcommand{\ee}{\end{equation}}
\newcommand{\bea}{\begin{eqnarray}\label}
\newcommand{\eea}{\end{eqnarray}}
\newcommand{\bas}{\begin{eqnarray*}}
\newcommand{\eas}{\end{eqnarray*}}
\newcommand{\bit}{\begin{itemize}}
\newcommand{\eit}{\end{itemize}}
\newcommand{\qed}{\hfill$\Box$ \vskip.2cm}
\newcommand{\nn}{\nonumber}
\newcommand{\R}{\mathbb{R}}
\newcommand{\N}{\mathbb{N}}
\newcommand{\pO}{\partial\Omega}

\newcommand{\eps}{\varepsilon}

\newcommand{\abs}{\\[5pt]}
\newcommand{\Abs}{\\[5mm]}

\newcommand{\io}{\int_\Omega}
\newcommand{\tm}{T_{max}}

\newcommand{\ophi}{\overline{\varphi}}
\newcommand{\uphi}{\underline{\varphi}}
\newcommand{\parab}{{\cal P}}
\newcommand{\qarab}{{\cal Q}}

\def\pa{\cdotp}
 
\def\va{\raise 2pt\hbox{,}}

\begin{document}
\title{A degenerate chemotaxis system with flux limitation:\\
 Maximally extended solutions and absence of gradient blow-up}
\author{
Nicola Bellomo\footnote{nicola.bellomo@polito.it}\\
{\small Department of Mathematics, Faculty  Sciences,}\\
 {\small King Abdulaziz University,  Jeddah, Saudi Arabia}\\
      {\small Politecnico di Torino,  10129 Torino, Italy}\\
\and
Michael Winkler\footnote{michael.winkler@math.uni-paderborn.de}\\
{\small Institut f\"ur Mathematik, Universit\"at Paderborn,}\\
{\small 33098 Paderborn, Germany} }
\date{}
\maketitle
\begin{abstract}
\noindent
  This paper aims at providing a first step toward a qualitative theory for
  a new class of chemotaxis models derived from the celebrated
  Keller-Segel system, with the main novelty being that diffusion is nonlinear with flux delimiter features.
  More precisely, as a prototypical representative of this class we study radially symmetric solutions
  of the parabolic-elliptic system
  \bas
	\left\{ \begin{array}{l}
\displaystyle u_t=\nabla \cdot \Big(\frac{u\nabla u}{\sqrt{u^2+|\nabla u|^2}}\Big) - \chi \, \nabla \cdot \Big(\frac{u\nabla v}{\sqrt{1+|\nabla v|^2}}\Big), \\[1mm]
	0=\Delta v - \mu + u,
	\end{array} \right.
  \eas
  under the initial condition $u|_{t=0}=u_0>0$ and no-flux boundary conditions
  in balls $\Omega\subset\R^n$, where $\chi>0$ and $\mu:=\frac{1}{|\Omega|} \io u_0$.\abs

  The main results assert the existence of a unique classical solution, extensible in time up to
  a maximal $\tm \in (0,\infty]$ which has the property that
  \bas
	\mbox{if} \quad \tm<\infty \quad \mbox{then} \quad
	\limsup_{t\nearrow\tm} \|u(\cdot,t)\|_{L^\infty(\Omega)}=\infty.
	\qquad \qquad (\star)
  \eas
  The proof therefor is mainly based on comparison methods,
  which firstly relate pointwise lower and upper bounds for the spatial gradient $u_r$
  to $L^\infty$ bounds for $u$ and to {\em upper bounds} for $z:=\frac{u_t}{u}$;
  secondly, another comparison argument involving nonlocal nonlinearities provides an appropriate control of
  $z_+$ in terms of bounds for $u$ and $|u_r|$,
  with suitably mild dependence on the latter.\abs
  As a consequence of ($\star$),
  by means of suitable a priori estimates it is moreover shown that the above solutions are global and bounded
  when either
  \bas
	n\ge 2 \ \mbox{ and } \chi<1,
	\qquad \mbox{or} \qquad
	n=1, \ \chi>0 \ \mbox{ and } m<m_c,
  \eas
  with $m_c:=\frac{1}{\sqrt{\chi^2-1}}$ if $\chi>1$ and $m_c:=\infty$ if $\chi\le 1$.\\
  That these conditions are essentially optimal will be shown in a forthcoming paper in which ($\star$)
  will be used to derive complementary results on the occurrence of solutions blowing up in finite time
  with respect to the norm of $u$ in $L^\infty(\Omega)$.\abs
\end{abstract}
\mysection{Introduction}
{\bf Keller-Segel systems with flux limitation.} \quad
The celebrated model by Keller and Segel \cite{Keller-Segel-70,Keller-Segel-71} was heuristically derived to model growth phenomena mediated by a chemoattractant, specifically the aggregation of dictyostelium discoideum due to an attractive chemical substance. The general structure of the model is as follows:
\begin{equation}
\label{KSO}
\left\{
\begin{array}{l}
  u_t  =
 \nabla \, \big(D_u(u,v) \nabla u - S(u,v) u \nabla v\big)+H_1(u,v), {} \\ [15pt]
  v_t   = D_v \Delta v + H_2(u,v),
\end{array} \right.
\end{equation}
where $u=u(x,t)$ denotes the cell (or organism) density at position $x$ and time $t$, and $v=v(x,t)$ is the
density of the chemoattractant. Here the function $S$ measures the chemotactic sensitivity,
the positive functions $D_u$ and $D_v$ represent the diffusivity of the cells and of the
chemoattractant, respectively, and $H_1$ and $H_2$ model source terms related to interactions. In a more general framework in which diffusions are not isotropic, $D_u$ and $S$ can be positive definite matrices.  See the survey by Hillen and Painter \cite{Hillen-Painter-09} for a review of modeling issues based on the classical approach of continuum mechanics closed by empirical models for the closure of conservation equations. The essay by Horstmann \cite{Horstmann-03} provides an additional source of information concerning modeling and applications in biology.
The recent survey \cite{BBYW} gives a review and qualitative analysis of a variety of mathematical problems and multiscale derivations of the
original model as well as of some recent developments such as the specific one treated in this paper. \abs
On the other hand, a natural question can be posed, namely if the use of parabolic models is consistent with the physics of the phenomena under consideration, or if, for instance, the use of hyperbolic models can be more appropriate. Or even within the approach by parabolic equations, if linear models are acceptable, while in the nonlinear case whether one should consider degenerate parabolic equations characterized by a finite propagation velocity. Intuitively, the answer is that phenomena with finite propagation velocity should be captured by an appropriate choice of nonlinear diffusion terms \cite{Vazquez-06}. A conceivable approach leads to consider functions $D_u(u,v)$ and $D_v(u,v)$ not only depending on $u$ and $v$, but also on their derivatives in space and time.
A recent study in this direction \cite{BBNS-10} has shown that macroscopic models can be obtained from the underlying description at the scale of cells delivered by suitable developments of kinetic theory methods. More in details, appropriate models of cell-cell interaction lead to macroscopic expressions for diffusion and cross-diffusion
with nonlinear limited flux terms of the type
\bas
	\nabla \cdot \bigg( D_u(u,v) \frac{ u\nabla u}{\sqrt{u^2 + \frac{\nu^2}{c^2}|\nabla u|^2}}\bigg)
	\qquad \hbox{and} \qquad
	\nabla \cdot \bigg(S(u,v) \frac{u\nabla v}{\sqrt{1 + |\nabla v|^2}} \bigg)\va
\eas
respectively,
with $\nu$ denoting the kinematic viscosity and $c$ the maximum speed of propagation, so that
in combination with an adequate equation for the evolution of the chemoattractant, this appraoch suggests to consider
models of type
\begin{equation}\label{Flim-KS}
\left\{
\begin{array}{l}
\displaystyle u_t =  \nabla \cdot \bigg(D_u (u,v)\frac{u \nabla u}{\sqrt{u^2 + \frac{\nu^2}{c^2}|\nabla u|^2}} -  S(u,v) \frac{u \nabla v}{\sqrt{1+|\nabla {v}|^2}}\bigg)
+H_1(u,v)\va \\[5mm]
\displaystyle v_t = D_v \Delta v + H_2(u,v),
\end{array} \right.
\end{equation}
as consistent modifications of the classical Keller-Segel system.
This idea, which is somehow related to the optimal transport framework \cite{Brenier-09}, can be motivated by a natural assumption of cell dynamics, where overcrowding is naturally avoided \cite{Burger-06}.
Furthermore, the introduction of this type of terms is founded in the assumption that
particles do not diffuse arbitrarily in space but, on the contrary, move through some privileged ways such as the border of cells.
In this new approach the non-physical diffusion is eliminated and the population moves with a finite speed of propagation, which is one of the intrinsic characteristics.
Indeed, the qualitative analysis of related systems with limited flux \cite{ACM-05,ACM-06} as well as some extensions to biological contexts (transport of morphogens) has been recently explored \cite{ACMS-12},
inter alia confirming the expected movement of fronts at finite speeds.\Abs
{\bf Boundedness vs.~blow-up.} \quad
In the framework of chemotaxis systems, however, a different qualitative aspect seems even more important, namely
the ability of the respective system to spontaneously generate structures.
In this regard, the classical Keller-Segel system, as obtained from (\ref{KSO}) on letting
$D_u\equiv D_v\equiv S\equiv 1$, $H\equiv 0$ and $K(u,v)=u-v$, is known to have the property that
some solutions reflect such aggregation processes even in the extreme mathematical sense of finite-time
blow-up of some solutions when either the spatial dimension $n$ satisfies $n\ge 3$ \cite{win_JMPA},
or when $n=2$ and the total mass of cells is suitably large \cite{herrero_velazquez, mizoguchi_win};
on the other hand, if either $n\ge 3$ and the initial data fulfill appropriate smallness conditions,
or $n=2$ and $\io u(\cdot,0)$ is small, or if $n=1$, then for various types of initial-boundary value problems,
global bounded solutions are known to exist
\cite{corrias_perthame, nagai_senba_yoshida, win_JDE2010, osaki_yagi, cao_optspace}.\abs
As for Keller-Segel-type models with flux limitations, the corresponding problem appears to be unsolved,
and it is the goal of the present paper to present a first step into a qualitative theory for such systems,
with a particular focus on the question whether solutions exist globally, as conjectured in \cite{BBNS-12},
or whether blow-up in finite time may occur for some initial data.
Specifically, we will consider the apparently most prototypical among the systems (\ref{Flim-KS}) in its
parabolic-elliptic simplification, as suggested in \cite{jaeger_luckhaus};
more precisely, we shall be concerned with the initial-boundary value problem
\be{0}
	\left\{ \begin{array}{l}
\displaystyle u_t=\nabla \cdot \bigg(\frac{u\nabla u}{\sqrt{u^2+|\nabla u|^2}}\Big) - \chi \, \nabla \cdot \Big(\frac{u\nabla v}{\sqrt{1+|\nabla v|^2}}\bigg)\va
	\qquad x\in \Omega, \ t>0, \\[5mm]
	0=\Delta v - \mu + u, \qquad x\in \Omega, \ t>0, \\[5mm]
	\displaystyle \bigg(\frac{u\nabla u}{\sqrt{u^2+|\nabla u|^2}} - \chi  \, \frac{u\nabla v}{\sqrt{1+|\nabla v|^2}} \bigg) \cdot \nu=0,
	\qquad x\in \partial\Omega, \ t>0, \\[5mm]
	u(x,0)=u_0(x), \qquad x\in\Omega,
 	\end{array} \right.
\ee
in a ball $\Omega=B_R(0)\subset \mathbb{R}^n$, $n\ge 1$, where
$\chi>0$ indicates the strength of chemotactic cross-diffusion.
In order to further simplify the analysis, we shall assume the initial data to satisfy
\be{init}
	u_0\in C^3(\bar\Omega)
	\quad \mbox{is radially symmetric and positive in $\bar\Omega$ with $\frac{\partial u_0}{\partial\nu}=0$ on } \pO,
\ee
so that the spatial average
\be{mu}
	\mu:=\frac{1}{|\Omega|} \io u_0(x)dx
\ee
is positive.\Abs
{\bf Main results.} \quad
In this framework, the first of our main results asserts local existence of a uniquely determined
classical solution.
In its most crucial part, however, the following theorem furthermore provides the extensibility criterion
(\ref{43.1}) which will be of great importance
both for deriving global existence in Theorem \ref{theo55} below,
as well as for characterizing the asymptotic behavior of non-global solutions
near their blow-up time \cite{bellomo_winkler2}.
%
%
\begin{theo}\label{theo43}
  Suppose that $u_0$ complies with (\ref{init}). Then there exist $\tm \in (0,\infty]$ and a uniquely determined
  pair $(u,v)$ of positive radially symmetric functions
  $u\in C^{2,1}(\bar\Omega\times [0,\tm))$ and $v\in C^{2,0}(\bar\Omega\times [0,\tm))$
  which solve (\ref{0}) classically in $\Omega\times (0,\tm)$, and which are such that
  \be{43.1}
	\mbox{if} \quad \tm<\infty \quad \mbox{then} \quad
	\limsup_{t\nearrow\tm} \|u(\cdot,t)\|_{L^\infty(\Omega)}=\infty.
  \ee
\end{theo}
In particular, (\ref{43.1}) rules out
the occurrence of any {\em gradient blow-up} phenomenon in the present framework;
it is thus impossible that $\nabla u$ becomes unbounded in finite time, whereas $u$ itself remains bounded.
In view of the complex evolution mechanism in (\ref{0}), inter alia involving doubly degenerate diffusion,	
this conclusion seems far from trivial;
indeed, various types of gradient-dependent nonlinearities and degeneracies are known to enforce unboundedness
of gradients for some solutions even in scalar reaction-diffusion equations
\cite{Angenent-96,Li-Souplet-10,Stinner-Winkler-08}.
Moreover, the additionally present cross-diffusive interaction apparently rules out the accessibility of (\ref{0})
to most of the techniques well-established in contexts of scalar parabolic equations with diffusion degeneracies
of related type, such as e.g.~the mean curvature flow equation and derivatives thereof, among others \cite{Evans-Spruck-91,Evans-Spruck-92}, or  \cite{Bertsch-DalPasso-92}. \abs
A natural next goal consists in identifying circumstances under which the above solutions are global.
Going in this direction,
the second of our main results provides conditions on the parameter $\chi$ in (\ref{0}) and, when $n=1$,
on the mass level $m$, which turn out to be sufficient not only for global extensibility,
but also for uniform boundedness of all solutions emanating from initial data $u_0$ with $\io u_0=m$.
\begin{theo}\label{theo55}
  Assume that $u_0$ satisfies (\ref{init}), and that either
  \be{55.111}
	n\ge 2
	\qquad \mbox{and} \qquad
	\chi<1,
  \ee
  or
  \be{55.112}
	n=1, \quad \chi>0
	\qquad \mbox{and} \qquad
	\io u_0<m_c,
  \ee
  where in the case $n=1$ we have set
  \be{mc}
	m_c:=\left\{ \begin{array}{ll}
	\frac{1}{\sqrt{\chi^2-1}} \qquad & \mbox{if } \chi>1, \\[3mm]
	+\infty & \mbox{if } \chi \le 1.
	\end{array} \right.
  \ee
  Then the problem (\ref{0}) possesses a unique global classical solution $(u,v)\in C^{2,1}(\bar\Omega\times [0,\infty)) \times C^{2,0}(\bar\Omega\times [0,\infty))$  which is radially symmetric and such that for some $C>0$ we have
  \be{55.1}
	\|u(\cdot,t)\|_{L^\infty(\Omega)} \le C
	\quad \mbox{and} \quad
	\|v(\cdot,t)\|_{L^\infty(\Omega)} \le C
	\quad \mbox{for all } t>0.
  \ee
\end{theo}
As a first and
immediate conclusion thereof, we underline that when $\chi<1$, in stark contrast to the original
Keller-Segel model, the system (\ref{0}) does
not exhibit any critical mass phenomenon, nor any phenomenon of critical sizes of initial data with
respect to global existence of solutions.
Let us secondly mention that the conditions (\ref{55.111}) and (\ref{55.112}), as identified above,
are in fact essentially optimal for the obtained conclusion:
Indeed, in \cite{bellomo_winkler2} the picture will in this respect be basically completed	
by showing that
if $\chi>1$ then in both cases $n=1$ with $m>m_c$, and $n\ge 2$,
some initial data can be constructed such that the corresponding solutions will blow up in finite time.

Together with the latter, our results thus indicate that in comparison to the original Keller-Segel system,
the occurrence of a critical mass phenomenon is shifted from the two-dimensional to the one-dimensional
setting, whereas in the case $n\ge 2$ we rather encounter a {\em critical sensitivity phenomenon}
in that the size of $\chi$ becomes the crucial quantity to determine whether or not blow-up may happen.\Abs
{\bf Main ideas. Excluding gradient blow-up.} \quad
In view of the doubly degenerate structure of the diffusion operator
$\nabla \cdot \Big(\frac{u\nabla u}{\sqrt{u^2+|\nabla u|^2}}\Big)$ in (\ref{0}),
standard theory yields local existence and extensibility as long as $u$ remains uniformly positive
and both $u$ and $\nabla u$ remain bounded (Lemma \ref{lem21}), where thanks to our positivity assumption
on $u_0$, a corresponding lower bound for $u$ can readily be obtained (Lemma \ref{lem23}).\abs
The crucial part in the derivation of Theorem \ref{theo43} will thus consist in ruling out
the possibility of gradient blow-up,
and in our approach toward this we will substantially make use of the radial symmetry of our solutions:
Based on 		
two different interpretations of the equation satisfied by $u_r$ as linear inhomogeneous parabolic
equations (Lemma 2.3), under the standing assumption that $u$ is non-global but remains bounded we will first obtain
a uniform lower bound for $u_r$ by a comparison argument (Lemma \ref{lem35}), and thereafter develop this
into a bound for $|u_r|$ in Section \ref{sect5}.\abs
The latter step itself will involve the quantity $z:=\frac{u_t}{u}$, as known to be of great importance
on various types of different nonlinear diffusion equations \cite{aronson, Vazquez-06}.
In the present context, we shall see that 	
$u_r$ can indeed be controlled in terms of the {\em positive part} $z_+$ of $z$ through an inequality
of the form
\be{est}
	\|u_r(\cdot,t)\|_{L^\infty((0,R))} \le C \cdot \Big(1+\|z_+\|_{L^\infty((0,R)\times (0,t))} \Big)
	\qquad \mbox{for all } t\in (0,\tm),
\ee
where $\tm\in (0,\infty)$ denotes the maximal existence time (Corollary \ref{cor40}).
This will be achieved by splitting the interval $(0,R)$ in two parts and
first performing a testing procedure to estimate $u_r$ in the corresponding inner region in certain weighted
Lebesgue spaces and taking limits appropriately (Lemma \ref{lem32}),
whereupon a comparison argument in the associated outer region will complete the proof of (\ref{est})
(Lemma \ref{lem37}).\abs
In order to complete the proof of Theorem \ref{theo43} by providing a suitable estimate for $z_+$,
we shall make use of the observation that $z$ satisfies the {\em one-sided} nonlocal parabolic inequality
\bas
	z_t(r,t) \le {\cal L} z 	
	+ d \cdot \Big( 1+\|z_+\|_{L^\infty((0,R)\times (0,t))} \Big)
\eas
with some $d>0$ and some homogeneous linear elliptic operator ${\cal L}$ (Lemma \ref{lem39}).
In fact, employing a maximum principle-type argument will show that this implies a pointwise upper bound
for $z$ (Lemma \ref{lem41}), which in conjunction with (\ref{est}) will prove Theorem \ref{theo43}.\abs
Thanks to the mild extensibility criterion (\ref{43.1}) thus gained, the proof of Theorem \ref{theo55}
thus actually reduces to the derivation of suitable a priori bounds for solutions with respect to
the norm of $u$ in $L^\infty(\Omega)$. This will be accomplished in the respective cases detailed in
Theorem \ref{theo55} by means of an essentially straightforward adaptation of the Moser-Alikakos iteration
technique to the present setting in Section \ref{sect6}.
\mysection{Preliminaries}		
\subsection{Local existence and a first extensibility criterion}
To begin with, let us suitably reduce (\ref{0}), locally in time, so as to become accessible to standard
existence theory.
We thereby obtain the following result on local existence of a smooth solution to (\ref{0}), extensible
as long as such a reduction is possible.
As a by-product, this procedure yields the first basic extensibility criterion (\ref{ext_crit}) the improvement
of which will be the main objective of the subsequent Sections \ref{sect3}-\ref{sect5}.
\begin{lem}\label{lem21}
  Suppose that $u_0$ satisfies (\ref{init}). Then there exist $\tm \in (0,\infty]$ and 		
  a uniquely determined
  pair $(u.v)$ of radially symmetric positive functions
  \bea{21.0}
	u \in C^{2,1}(\bar\Omega\times [0,\tm)), \qquad v\in C^{2,0}(\bar\Omega\times [0,\tm)),
  \eea
  which solve (\ref{0}) classically in $\Omega\times (0,\tm)$, and which are such that
  \be{ext_crit}
	\mbox{if $\tm<\infty$ \quad then either \quad }
	\liminf_{t\nearrow\tm} \inf_{x\in\Omega} u(x,t)=0 \mbox{ \quad or \quad }
	\limsup_{t\nearrow \tm} \|u(\cdot,t)\|_{W^{1,\infty}(\Omega)}=\infty.
  \ee
\end{lem}
\proof
  We let
  \be{21.1}
	\eps:=\min \bigg\{ \frac{1}{2} \inf_{x\in\Omega} u_0(x) \, , \, \frac{1}{2\|u_0\|_{L^\infty(\Omega)}} \, , \,
	\frac{1}{2\|\nabla u_0\|_{L^\infty(\Omega)}} \bigg\}
  \ee
  and take cut-off functions $\psi_\eps \in C^\infty(\R)$ and $\phi_\eps\in C^\infty(\R)$ satisfying
  \bas
	\frac{\eps}{2} \le \psi_\eps(s) \le \frac{2}{\eps}
	\quad \mbox{for all } s\in\R
	\qquad \mbox{and} \qquad
	\psi_\eps(s)=s
	\quad \mbox{for all } s\in \Big( \eps \, , \, \frac{1}{\eps} \Big)\va
  \eas
  as well as
  \bas
	\phi_\eps(s)\le \frac{2}{\eps}
	\quad \mbox{for all } s\in\R
	\qquad \mbox{and} \qquad
	\phi_\eps(s)=s
	\quad \mbox{for all } s\le \frac{1}{\eps}\pa
  \eas

  Then
  \bas
	a_\eps(s,p):=\frac{\psi_\eps(s)}{\sqrt{\psi_\eps^2(s)+\phi_\eps^2(|p|)}},
	\qquad s\in\R, \ p\in\R^n,
  \eas
  defines a function $a_\eps \in C^\infty(\R\times \R^n)$ fulfilling
  \bas
	a_\eps(s,p) \le \frac{\psi_\eps(s)}{\sqrt{\psi_\eps^2(s)}}=1
	\quad \mbox{for all $s\in\R$ and } p\in\R^n
  \eas
  and
  \bas
	a_\eps(s,p) \ge \frac{\frac{\eps}{2}}{\sqrt{(\frac{2}{\eps})^2 + (\frac{2}{\eps})^2}} = \frac{\eps^2}{4\sqrt{2}}
	\quad \mbox{for all $s\in\R$ and } p\in\R^n.
  \eas
  We can therefore adapt a fixed point argument which is well-established in the existence theory of parabolic-elliptic
  chemotaxis systems (cf.~\cite{Djie-Winkler-10}  or \cite{Fujie-Winkler-Yokota}, for instance) to find $T_\eps>0$ such that the problem
  \bas
	\left\{ \begin{array}{ll}
\displaystyle	u_t=\nabla \cdot \Big( \frac{\psi_\eps(u)\nabla u}{\sqrt{\psi_\eps^2(u) + \phi_\eps^2(|\nabla u|)}}\Big)
	- \chi \nabla \cdot \Big( \frac{u\nabla v}{\sqrt{1+|\nabla v|^2}}\Big)\va \qquad x\in\Omega, \ t\in (0,T_\eps), \\[5mm]
	0=\Delta v - \mu + u, \qquad x\in\Omega, \ t\in (0,T_\eps), \\[5mm]
\displaystyle	\frac{\partial u}{\partial\nu}=\frac{\partial v}{\partial\nu}=0, 	\qquad x\in\partial\Omega, \ t\in (0,T_\eps), \\[5mm]
	u(x,0)=u_0(x),
	\qquad x\in\Omega,
	\end{array} \right.
  \eas
  possesses a unique classical solution $(u_\eps,v_\eps)$ such that $u_\eps \in C^{2,1}(\bar\Omega \times [0,T_\eps))$
  and $v_\eps \in C^{2,0}(\bar\Omega\times [0,T_\eps))$, and such that both $u_\eps$ and $v_\eps$ are radially symmetric and
  positive. Furthermore, since
  $2\eps \le u_0 \le \frac{1}{2\eps}$ and $|\nabla u_0| \le \frac{1}{2\eps}$ in $\Omega$ according to our choice of $\eps$,
  by continuity of $u_\eps$ and $\nabla u_\eps$ in $\bar\Omega \times [0,T_\eps)$ we can find  $\tilde T_\eps \in (0,T_\eps)$
  such that
  \bas
	\eps \le u_\eps \le \frac{1}{\eps}
	\quad \mbox{and} \quad
	|\nabla u_\eps| \le \frac{1}{\eps}
	\qquad \mbox{in } \Omega \times (0,\tilde T_\eps).
  \eas
  In particular, this implies that $\psi_\eps(u_\eps)=u_\eps$ and $\phi_\eps(|\nabla u_\eps|)=|\nabla u_\eps|$ in
  $\Omega\times (0,\tilde T_\eps)$, and that thus $a_\eps(u_\eps)=\frac{u_\eps}{\sqrt{u_\eps^2 + |\nabla u_\eps|^2}}$
  in this region, meaning that $(u_\eps,v_\eps)$ actually solves the original problem (\ref{0}) in
  $\Omega\times (0,\tilde T_\eps)$.

  Finally, in view of the dependence of $\eps$ on $u_0$ as expressed in (\ref{21.1}),
  a standard extensibility argument yields that
  the above solution can be continued so as to exist up to some maximal time $\tm \le \infty$ in such a way that
  (\ref{ext_crit}) is valid.
\qed
\subsection{Radial solutions}
Since all our solutions are radially symmetric, whenever this appears convenient me may without any danger of confusion
utilize the notation $u(r,t)$ and $v(r,t)$
instead of $u(x,t)$ and $v(x,t)$, respectively, where $r=|x|\in (0,R)$.\abs
In this particular radial setting, $u$ actually fulfills a favorable parabolic equation specified in the following lemma.
\begin{lem}\label{lem36}
  Assume (\ref{init}). Then the solution of (\ref{0}) satisfies
  \bea{36.1}
	u_t &=& \frac{u^3 u_{rr}}{\sqrt{u^2+u_r^2}^3} + \frac{u_r^4}{\sqrt{u^2+u_r^2}^3}
		+ \frac{n-1}{r} \cdot \frac{uu_r}{\sqrt{u^2+u_r^2}} \nn\\
	& & - \chi \frac{u_r v_r}{\sqrt{1+v_r^2}} - \chi \frac{u(\mu-u)}{\sqrt{1+v_r^2}^3}
	- \chi \cdot \frac{n-1}{r} \cdot \frac{uv_r^3}{\sqrt{1+v_r^2}^3},
  \eea
  for all $r\in (0,R)$ and $t\in (0,\tm)$.
\end{lem}
\proof
  We differentiate on the right-hand side of the first equation in (\ref{0}) to obtain
  \bea{36.2}
	u_t &=& \frac{1}{r^{n-1}} \cdot \bigg( r^{n-1} \frac{uu_r}{\sqrt{u^2+u_r^2}}\bigg)_r
	- \frac{\chi}{r^{n-1}}  \cdot \bigg( r^{n-1} \frac{uv_r}{\sqrt{1+v_r^2}}\bigg)_r \nn\\
	&=& \frac{uu_{rr}}{\sqrt{u^2+u_r^2}} + \frac{u_r^2}{\sqrt{u^2+u_r^2}}
	-\frac{1}{2} \cdot \frac{uu_r(2uu_r+2u_r u_{rr})}{\sqrt{u^2+u_r^2}^3}
	+ \frac{n-1}{r} \cdot \frac{uu_r}{\sqrt{u^2+u_r^2}} \nn\\
	& & - \chi \frac{uv_{rr}}{\sqrt{1+v_r^2}} - \chi \frac{u_r v_r}{\sqrt{1+v_r^2}}
	+ \frac{1}{2} \cdot \chi \frac{uv_r \cdot 2v_r v_{rr}}{\sqrt{1+v_r^2}^3}
	-\chi  \cdot \frac{n-1}{r} \cdot \frac{uv_r}{\sqrt{1+v_r^2}}
  \eea
  for all $r\in (0,R)$ and $t\in (0,\tm)$.
  Here we can rearrange
  \bas
	& & \hspace*{-30mm}
	\frac{uu_{rr}}{\sqrt{u^2+u_r^2}} + \frac{u_r^2}{\sqrt{u^2+u_r^2}}
	-\frac{1}{2} \cdot \frac{uu_r(2uu_r+2u_r u_{rr})}{\sqrt{u^2+u_r^2}^3} \\
	&=& \frac{uu_{rr}}{\sqrt{u^2+u_r^2}^3} \cdot \Big\{ (u^2+u_r^2) - u_r^2 \Big\}
	+ \frac{u_r^2}{\sqrt{u^2+u_r^2}^3} \cdot \Big\{ (u^2+u_r^2)-u^2 \Big\} \\
	&=& \frac{u^3 u_{rr}}{\sqrt{u^2+u_r^2}^3} + \frac{u_r^4}{\sqrt{u^2+u_r^2}^3}
  \eas
  and, similarly,
  \bas
	& & \hspace*{-30mm}
	- \chi \frac{uv_{rr}}{\sqrt{1+v_r^2}}
	+ \frac{1}{2} \cdot \chi \frac{uv_r \cdot 2v_r v_{rr}}{\sqrt{1+v_r^2}^3}
	-\chi  \cdot \frac{n-1}{r} \cdot \frac{uv_r}{\sqrt{1+v_r^2}} \\
	&=& - \chi \frac{uv_{rr}}{\sqrt{1+v_r^2}^3} \cdot \Big\{ (1+v_r^2)-v_r^2 \Big\}
	- \chi \cdot \frac{n-1}{r} \cdot \frac{uv_r}{\sqrt{1+v_r^2}^3} \cdot (1+v_r^2) \\
	&=& - \chi \frac{u(v_{rr}+\frac{n-1}{r}v_r)}{\sqrt{1+v_r^2}^3}
	- \chi \cdot \frac{n-1}{r} \cdot \frac{uv_r^3}{\sqrt{1+v_r^2}^3}
  \eas
  for $r\in (0,R)$ and $t\in (0,\tm)$. Since $v_{rr}+\frac{n-1}{r}v_r=\mu-u$ by (\ref{0}), the identity (\ref{36.1})
  thus results from (\ref{36.2}).
\qed
We next differentiate (\ref{36.1}) to obtain a corresponding equation for $u_r$.
Here suitable arrangements will lead to the two alternative interpretations (\ref{35.01}) and (\ref{35.02})
thereof as linear inhomogeneous parabolic equations.
The first of these will be used to establish an estimate {\em from below} for $u_r$ in Lemma \ref{lem34} by a straightforward
comparison argument, whereas upon some more involved preparations, on the basis of the latter we will
apply another comparison procedure to derive a certain {\em upper} bound for $u_r$ in Lemma \ref{lem37}.
\begin{lem}\label{lem35}
  Assume (\ref{init}). Then
  \bea{35.1}
	u_{rt}
	&=& \frac{u^3 u_{rrr}}{\sqrt{u^2+u_r^2}^3}
	+ 3\frac{u^2 u_r^3 u_{rr}}{\sqrt{u^2+u_r^2}^5}
	- 3\frac{u^3 u_r u_{rr}^2}{\sqrt{u^2+u_r^2}^5} \nn\\
	& & + 4\frac{u^2 u_r^3 u_{rr}}{\sqrt{u^2+u_r^2}^5}
	+ \frac{u_r^5 u_{rr}}{\sqrt{u^2+u_r^2}^5}
	-3\frac{uu_r^5}{\sqrt{u^2+u_r^2}^5} \nn\\
	& & - \frac{n-1}{r^2} \cdot \frac{uu_r}{\sqrt{u^2+u_r^2}}
	+ \frac{n-1}{r} \cdot \frac{u^3 u_{rr}}{\sqrt{u^2+u_r^2}^3}
	+ \frac{n-1}{r} \frac{u_r^4}{\sqrt{u^2+u_r^2}^3} \nn\\
	& & - \chi\mu \frac{u_r}{\sqrt{1+v_r^2}^3}
	+ 2\chi \frac{uu_r}{\sqrt{1+v_r^2}^3}
	+ 3\chi\mu \frac{uv_r v_{rr}}{\sqrt{1+v_r^2}^5}
	- 3\chi \frac{u^2 v_r v_{rr}}{\sqrt{1+v_r^2}^5} \nn\\
	& & -\chi \frac{u_{rr} v_r}{\sqrt{1+v_r^2}}
	-\chi \frac{u_r v_{rr}}{\sqrt{1+v_r^2}}
	+ \chi \frac{u_r v_r^2 v_{rr}}{\sqrt{1+v_r^2}^3} \nn\\
	& & + \chi \cdot \frac{n-1}{r^2} \cdot \frac{uv_r^3}{\sqrt{1+v_r^2}^3}
	- \chi \cdot \frac{n-1}{r} \cdot \frac{u_r v_r^3}{\sqrt{1+v_r^2}^3}
	-3\chi \cdot \frac{n-1}{r} \cdot \frac{uv_r^2 v_{rr}}{\sqrt{1+v_r^2}^5}
  \eea
  for all $r\in (0,R)$ and $t\in (0,\tm)$. In particular,
  \be{35.01}
	(\parab u_r)(r,t)=0
	\quad \mbox{for all $r\in (0,R)$ and } t\in (0,\tm),
  \ee
  where the inhomogeneous linear parabolic operator $\parab$ is defined by
  \be{parab}
	(\parab \varphi)(r,t)
	:= \varphi_t - A_1(r,t)\varphi_{rr} - A_2(r,t) \varphi_r - A_3(r,t)\varphi - A_4(r,t),
	\qquad r\in (0,R), \ t\in (0,\tm),
  \ee
  with
  \bea{parab1}
	A_1(r,t) &:=& \frac{u^3}{\sqrt{u^2+u_r^2}^3}, \nn\\[1mm]
	A_2(r,t)&:=& 3\frac{u^2 u_r^3}{\sqrt{u^2+u_r^2}^5}
	- 3\frac{u^3 u_r u_{rr}}{\sqrt{u^2+u_r^2}^5}
	+ 4\frac{u^2 u_r^3}{\sqrt{u^2+u_r^2}^5}
	+ \frac{u_r^5}{\sqrt{u^2+u_r^2}^5}
	+ \frac{n-1}{r} \cdot \frac{u^3}{\sqrt{u^2+u_r^2}^3} \nn\\
	& & - \chi \frac{v_r}{\sqrt{1+v_r^2}}, \nn\\[1mm]
	A_3(r,t) &:=& -3\frac{uu_r^4}{\sqrt{u^2+u_r^2}^5}
	- \frac{n-1}{r^2} \frac{u}{\sqrt{u^2+u_r^2}} \nn\\
	& & -\chi\mu \frac{1}{\sqrt{1+v_r^2}^3}
	+2\chi \frac{u}{\sqrt{1+v_r^2}^3}
	- \chi \frac{v_{rr}}{\sqrt{1+v_r^2}}
	+ \chi \frac{v_r^2 v_{rr}}{\sqrt{1+v_r^2}^3}
	- \chi \cdot \frac{n-1}{r} \cdot \frac{v_r^3}{\sqrt{1+v_r^2}^3}
	\qquad \mbox{and} \nn\\[1mm]
	A_4(r,t) &:=& \frac{n-1}{r} \cdot \frac{u_r^4}{\sqrt{u^2+u_r^2}^3} \nn\\
	& & + 3\chi\mu \frac{u v_r v_{rr}}{\sqrt{1+v_r^2}^5}
	-3\chi \frac{u^2 v_r v_{rr}}{\sqrt{1+v_r^2}^5}
	+ \chi \cdot \frac{n-1}{r^2} \cdot \frac{u v_r^3}{\sqrt{1+v_r^2}^3}
	-3\chi \cdot \frac{n-1}{r} \cdot \frac{uv_r^2 v_{rr}}{\sqrt{1+v_r^2}^5}
  \eea
  for $r\in (0,R)$ and $t\in (0,\tm)$.
  Likewise,
  \be{35.02}
	(\qarab u_r)(r,t)=0
	\qquad \mbox{for all $r\in (0,R)$ and } t\in (0,\tm),
  \ee
  with $\qarab$ given by
  \be{qarab}
	(\qarab \varphi)(r,t):=\varphi_t - A_1(r,t) \varphi_{rr} - A_2(r,t)\varphi_r
	-\tilde A_3(r,t) \varphi - \tilde A_4(r,t),
	\qquad r\in (0,R), \ t\in (0,\tm),
  \ee
  where
  \bea{qarab1}
	\tilde A_3(r,t)
	&:=& \frac{n-1}{r} \cdot \frac{u_r^3}{\sqrt{u^2+u_r^2}^3} \nn\\
	& & - \chi \mu \frac{1}{\sqrt{1+v_r^2}^3}
	+ 2\chi \frac{u}{\sqrt{1+v_r^2}^3}
	-\chi \frac{v_{rr}}{\sqrt{1+v_r^2}}
	+ \chi \frac{v_r^2 v_{rr}}{\sqrt{1+v_r^2}^3}
	- \chi \cdot \frac{n-1}{r} \cdot \frac{v_r^3}{\sqrt{1+v_r^2}^3}
	\qquad \mbox{and} \nn\\[1mm]
	\tilde A_4(r,t)
	&:=& -3\frac{uu_r^5}{\sqrt{u^2+u_r^2}^5}
	- \frac{n-1}{r^2} \cdot \frac{uu_r}{\sqrt{u^2+u_r^2}} \nn\\
	& & + 3\chi\mu \frac{uv_r v_{rr}}{\sqrt{1+v_r^2}^5}
	-3\chi \frac{u^2 v_r v_{rr}}{\sqrt{1+v_r^2}^5}
	+ \chi \cdot \frac{n-1}{r} \cdot \frac{uv_r^3}{\sqrt{1+v_r^2}^3}
	- 3\chi \cdot \frac{n-1}{r} \cdot \frac{uv_r^2 v_{rr}}{\sqrt{1+v_r^2}^5}
  \eea
  for $r\in (0,R)$ and $t\in (0,\tm)$.
\end{lem}
\proof
  Differentiation of (\ref{36.1}) with respect to $r$ yields
  \bea{35.2}
	u_{rt}
	&=& \frac{u^3 u_{rrr}}{\sqrt{u^2+u_r^2}^3}
	+ 3\cdot \frac{u^2 u_r u_{rr}}{\sqrt{u^2+u_r^2}^3}
	- \frac{3}{2} \cdot \frac{u^3 u_{rr} \cdot (2uu_r+2u_r u_{rr})}{\sqrt{u^2+u_r^2}^5} \nn\\
	& & + 4\cdot \frac{u_r^3 u_{rr}}{\sqrt{u^2+u_r^2}^3}
	-\frac{3}{2} \cdot \frac{u_r^4 \cdot (2uu_r+2u_r u_{rr})}{\sqrt{u^2+u_r^2}^5} \nn\\
	& & -\frac{n-1}{r^2} \cdot \frac{uu_r}{\sqrt{u^2+u_r^2}}
	+ \frac{n-1}{r} \cdot \frac{uu_{rr}}{\sqrt{u^2+u_r^2}} \nn\\
	& & +\frac{n-1}{r} \cdot \frac{u_r^2}{\sqrt{u^2+u_r^2}}
	-\frac{1}{2} \cdot \frac{n-1}{r} \cdot \frac{uu_r \cdot (2uu_r + 2u_r u_{rr})}{\sqrt{u^2+u_r^2}^3} \nn\\
	& & -\chi\mu \cdot \frac{u_r}{\sqrt{1+v_r^2}^3}
	+ 2\chi \cdot \frac{uu_r}{\sqrt{1+v_r^2}^3}
	+\frac{3}{2} \chi  \cdot\frac{u(\mu-u) \cdot 2v_r v_{rr}}{\sqrt{1+v_r^2}^5} \nn\\
	& & -\chi \cdot \frac{u_{rr} v_r}{\sqrt{1+v_r^2}}
	-\chi \cdot \frac{u_r v_{rr}}{\sqrt{1+v_r^2}}
	+ \frac{1}{2} \chi \cdot \frac{u_r v_r \cdot 2v_r v_{rr}}{\sqrt{1+v_r^2}^3} \nn\\
	& & +\chi \cdot \frac{n-1}{r^2} \cdot \frac{uv_r^3}{\sqrt{1+v_r^2}^3}
	-\chi \cdot \frac{n-1}{r} \cdot \frac{u_r v_r^3}{\sqrt{1+v_r^2}^3} \nn\\
	& & -3\chi \cdot \frac{n-1}{r} \cdot \frac{u v_r^2 v_{rr}}{\sqrt{1+v_r^2}^3}
	+\frac{3}{2}\chi \cdot \frac{n-1}{r} \cdot \frac{uv_r^3 \cdot 2v_r v_{rr}}{\sqrt{1+v_r^2}^5}
  \eea
  for $r\in (0,R)$ and $t\in (0,\tm)$, where
  \bas
	3 \cdot\frac{u^2 u_r u_{rr}}{\sqrt{u^2+u_r^2}^3}
	- \frac{3}{2} \cdot \frac{u^3 u_{rr} \cdot (2uu_r+2u_r u_{rr})}{\sqrt{u^2+u_r^2}^5}
	&=& 3 \cdot\frac{u^2 u_r u_{rr}}{\sqrt{u^2+u_r^2}^5} \cdot
	\Big\{ (u^2+u_r^2)-u^2\Big\}
	-3 \cdot\frac{u^3 u_r u_{rr}^2}{\sqrt{u^2+u_r^2}^5} \\
	&=& 3 \cdot\frac{u^2 u_r^3 u_{rr}}{\sqrt{u^2+u_r^2}^5}
	-3 \cdot\frac{u^3 u_r u_{rr}^2}{\sqrt{u^2+u_r^2}^5}
  \eas
  and
  \bas
	4 \cdot\frac{u_r^3 u_{rr}}{\sqrt{u^2+u_r^2}^3}
	-\frac{3}{2} \cdot \frac{u_r^4 \cdot (2uu_r+2u_r u_{rr})}{\sqrt{u^2+u_r^2}^5}
	&=& \frac{u_r^3 u_{rr}}{\sqrt{u^2+u_r^2}^5} \cdot \Big\{ 4(u^2+u_r^2)-3u_r^2\Big\}
	- 3 \cdot\frac{uu_r^5}{\sqrt{u^2+u_r^2}^5} \\
	&=& 4 \cdot\frac{u^2 u_r^3 u_{rr}}{\sqrt{u^2+u_r^2}^5}
	+ \frac{u_r^5 u_{rr}}{\sqrt{u^2+u_r^2}^5}
	- 3 \cdot\frac{uu_r^5}{\sqrt{u^2+u_r^2}^5}
  \eas
  as well as
  \bas
	& & \hspace*{-20mm}
	\frac{n-1}{r} \cdot \frac{uu_{rr}}{\sqrt{u^2+u_r^2}}
	+\frac{n-1}{r} \cdot \frac{u_r^2}{\sqrt{u^2+u_r^2}}
	-\frac{1}{2} \cdot \frac{n-1}{r} \cdot \frac{uu_r \cdot (2uu_r + 2u_r u_{rr})}{\sqrt{u^2+u_r^2}^3} \\
	&=& \frac{n-1}{r} \cdot \frac{uu_{rr}}{\sqrt{u^2+u_r^2}^3} \cdot \Big\{ (u^2+u_r^2)-u_r^2 \Big\}
	+ \frac{n-1}{r} \cdot \frac{u_r^2}{\sqrt{u^2+u_r^2}^3} \cdot \Big\{(u^2+u_r^2)-u^2\Big\} \\
	&=& \frac{n-1}{r} \cdot \frac{u^3 u_{rr}}{\sqrt{u^2+u_r^2}^3} 	
	+\frac{n-1}{r} \cdot \frac{u_r^4}{\sqrt{u^2+u_r^2}^3}
  \eas
  for $r\in (0,R)$ and $t\in (0,\tm)$.
  Finally simplifying the last two summands in (\ref{35.2}) according to
  \bas
	& & \hspace*{-30mm}
	-3\chi \cdot \frac{n-1}{r} \cdot \frac{u v_r^2 v_{rr}}{\sqrt{1+v_r^2}^3}
	+\frac{3}{2}\chi \cdot \frac{n-1}{r} \cdot \frac{uv_r^3 \cdot 2v_r v_{rr}}{\sqrt{1+v_r^2}^5} \\
	&=& -3\chi \cdot \frac{n-1}{r} \cdot \frac{uv_r^2 v_{rr}}{\sqrt{1+v_r^2}^5} \cdot \Big\{(1+v_r^2) - v_r^2 \Big\} \\
	&=& -3\chi \cdot \frac{n-1}{r} \cdot \frac{uv_r^2 v_{rr}}{\sqrt{1+v_r^2}^5} \va
  \eas
  from (\ref{35.2}) we easily obtain (\ref{35.1}), and thus also (\ref{35.01}) and (\ref{35.02}).
\qed
Thanks to the favorable structure of the equation for $v$ in (\ref{0}), this second solution component
can be expressed explicitly in terms of $u$. This leads to the following observations which will frequently
be referred to throughout the sequel.
\begin{lem}\label{lem33}
  Assume (\ref{init}). Then
  \be{33.1}
	v_r(r,t)=\frac{\mu r}{n} - r^{1-n} \cdot \int_0^r \rho^{n-1} u(\rho,t)d\rho
	\quad \mbox{for all $r\in (0,R)$ and } t\in (0,\tm)
  \ee
  and
  \be{33.2}
	v_{rr}(r,t)=\frac{\mu}{n} - u + \frac{n-1}{r^n} \cdot \int_0^r \rho^{n-1} u(\rho,t)d\rho
	\quad \mbox{for all $r\in (0,R)$ and } t\in (0,\tm).
  \ee
  Moreover we have
  \be{33.3}
	v_{rt}= - \frac{uu_r}{\sqrt{u^2+u_r^2}} + \chi \cdot \frac{uv_r}{\sqrt{1+v_r^2}}
	\quad \mbox{in } (0,R)\times (0,\tm).
  \ee
\end{lem}
\proof
  Since by the second equation in (\ref{0}) we have
  \bas
	(r^{n-1}v_r)_r=r^{n-1}(\mu-u)
	\quad \mbox{for all $r\in (0,R)$ and } t\in (0,\tm),
  \eas
  the identity (\ref{33.1}) easily results by integration, whereupon a differentiation of (\ref{33.1}) with respect
  to $r$ yields (\ref{33.2}).\\
  Next we differentiate (\ref{33.1}) with respect to $t$ and use the first equation in (\ref{0}) to see that
  \bas
	v_{rt}(r,t)
	&=& - \frac{1}{r^{n-1}} \cdot \int_0^r \rho^{n-1} u_t(\rho,t)d\rho \\
	&=& - \frac{1}{r^{n-1}} \cdot \int_0^r \rho^{n-1} \cdot
	\bigg\{ \frac{1}{\rho^{n-1}} \Big( \rho^{n-1} \frac{uu_r}{\sqrt{u^2+u_r^2}}
		-\chi \rho^{n-1} \frac{uv_r}{\sqrt{1+v_r^2}} \Big) \bigg\}_r (\rho,t) d\rho \\
	&=& - \frac{1}{r^{n-1}} \cdot \Big\{ r^{n-1} \frac{uu_r}{\sqrt{u^2+u_r^2}}
		-\chi r^{n-1} \frac{uv_r}{\sqrt{1+v_r^2}}\Big\}
  \eas
  for all $r\in (0,R)$ and $t\in (0,\tm)$, which shows (\ref{33.3}).
\qed
Let us note some pointwise estimates resulting from Lemma \ref{lem333} in a straightforward manner.
\begin{lem}\label{lem333}
  Let (\ref{init}) hold. Then for each $t\in (0,\tm)$ and any $r\in (0,R)$ we have
  \be{333.0}
	-\frac{\mu R^n}{n} \cdot r^{1-n} \le v_r(r,t) \le \frac{\mu}{n} \cdot r
  \ee
  and
  \be{333.1}
	|v_r(r,t)| \le \frac{\|u(\cdot,t)\|_{L^\infty((0,R))}}{n} \cdot r
  \ee
  as well as
  \be{333.2}
	|v_{rr}(r,t)| \le \|u(\cdot,t)\|_{L^\infty((0,R))}.
  \ee
\end{lem}
\proof
  Fixing $t\in (0,\tm)$ and writing $M:=\|u(\cdot,t)\|_{L^\infty((0,R))}$, we clearly have $\mu \le M$, so that since from
  Lemma \ref{lem33} we know that
  \bas
	v_r(r,t) \le \frac{\mu}{n} \cdot r
	\quad \mbox{for all } r\in (0,R)
  \eas
  and
  \bas
	v_r(r,t) \ge - \frac{1}{r^{n-1}} \cdot \int_0^r \rho^{n-1} \cdot M d\rho
	= - \frac{Mr}{n}
	\quad \mbox{for all} r\in (0,R),
  \eas
  both the right inequality in (\ref{333.0}) as well as (\ref{333.1}) are immediate.
  Similarly, using (\ref{33.2}) we can estimate
  \bas
	v_{rr}(r,t) &\le& \frac{\mu}{n} + \frac{n-1}{r^n} \cdot \int_0^r \rho^{n-1} \cdot M d\rho = \frac{\mu}{n} + \frac{(n-1)M}{n} \\
	&\le& M
	\quad \mbox{for all } r\in (0,R)
  \eas
  and
  \bas
	v_{rr}(r,t) \ge -M
	\quad \mbox{for all } r\in (0,R),
  \eas
  which yields (\ref{333.2}).
  Finally, to derive the left inequality in (\ref{333.0}) we observe that $\int_0^r \rho^{n-1} u(\rho,t) d\rho
  \le \frac{m}{\omega_n}$ for all $r\in (0,R)$
  and recall that $\frac{m}{\omega_n}=\frac{\mu R^n}{n}$ by (\ref{mu}) to obtain from (\ref{33.1}) that
  \bas
	v_r(r,t) \ge - r^{1-n} \int_0^r \rho^{n-1} u(\rho,t) d\rho
	\ge - r^{1-n} \cdot \frac{m}{\omega_n}
	= - r^{1-n} \cdot \frac{\mu R^n}{n}
	\quad \mbox{for all } r\in (0,R).
  \eas
  This completes the proof.
\qed

\mysection{A pointwise estimate from below for $u$}\label{sect3}
In order to show that (\ref{ext_crit}) actually reduces to (\ref{43.1}), let us first rule out the occurrence
of the first alternative in (\ref{ext_crit}).
In proving this, we shall make use of the following elementary inequality.
\begin{lem}\label{lem24}
  We have
  \bas
	\frac{\xi}{\sqrt{1+\xi}^3} \le \frac{2}{3\sqrt{3}}
	\qquad \mbox{for all } \xi\ge 0.
  \eas
\end{lem}
\proof
  Since $\varphi(\xi):=\frac{\xi}{\sqrt{1+\xi}^3}, \ \xi\ge 0$, satisfies $\varphi(0)=0, \varphi(\xi)\to 0$ as
  $\xi\to\infty$ and
  $\varphi'(\xi)=(1+\xi)^{-\frac{5}{2}} \cdot (1-\frac{\xi}{2})$ for all $\xi>0$, it follows that
  $\varphi(\xi) \le \varphi(2) = \frac{2}{\sqrt{3}^3}$ for all $\xi\ge 0$.
\qed
By means of a comparison argument applied to (\ref{36.1}),
we can now in fact exclude that solutions attain zeros within finite time.
\begin{lem}\label{lem23}
  If (\ref{init}) holds, then
  \be{23.1}
	u(r,t) \ge \Big( \inf_{r\in (0,R)} u_0(r) \Big) \cdot e^{-\kappa t}
	\qquad \mbox{for all $r\in (0,R)$ and } t\in (0,\tm),
  \ee
  where
  \be{23.2}
	\kappa:=\chi \mu + \frac{2(n-1)\chi\mu}{3\sqrt{3}n}.
  \ee
\end{lem}
\proof
  We rewrite (\ref{36.1}) in the form
  \be{23.3}
	u_t= a_1(r,t) u_{rr} + a_{21}(r,t)u_r + \frac{a_{22}(r,t)}{r} \cdot u_r
	-\chi \cdot \frac{u(\mu-u)}{\sqrt{1+v_r^2}^3}
	-\frac{n-1}{r} \cdot \chi \frac{uv_r^3}{\sqrt{1+v_r^2}^3}
  \ee
  for all $r\in (0,R)$ and $t\in (0,\tm)$, where
  \bas
	a_1(r,t):=\frac{u^3}{\sqrt{u^2+u_r^2}^3}
  \eas
  and
  \bas
	a_{21}(r,t):=\frac{u_r^3}{\sqrt{u^2+u_r^2}^3} - \chi \cdot \frac{v_r}{\sqrt{1+v_r^2}}
  \eas
  as well as
  \bas
	a_{22}(r,t):=(n-1)\cdot \frac{u}{\sqrt{u^2+u_r^2}}
  \eas
  define continuous functions in $[0,R] \times (0,\tm)$.
  In (\ref{23.3}), we can estimate
  \bas
	- \chi \cdot\frac{u(\mu-u)}{\sqrt{1+v_r^2}^3}
	\ge - \chi\mu \cdot \frac{u}{\sqrt{1+v_r^2}^3}
	\ge - \chi\mu u
	\qquad \mbox{for all $r\in (0,R)$ and } t\in (0,\tm),
  \eas
  and in order to control the last term in (\ref{23.3}) we
  use the one-sided inequality $v_r \le \frac{\mu r}{n}$ provided by Lemma \ref{lem333}, which in conjunction with
  Lemma \ref{lem24} entails that
  \bas
	- \frac{n-1}{r} \cdot \chi \frac{uv_r^3}{\sqrt{1+v_r^2}^3}
	&=& - (n-1)\chi \cdot \frac{v_r^2}{\sqrt{1+v_r^2}^3} \cdot \frac{v_r}{r} \cdot u \\
	&\ge& - (n-1)\chi \cdot \frac{2}{3\sqrt{3}} \cdot \frac{\mu}{n} \cdot u
	\qquad \mbox{for all $r\in (0,R)$ and } t\in (0,\tm).
  \eas
  Accordingly, from (\ref{23.3}) we infer that with $\kappa$ as in (\ref{23.2}) we have
  \bas
	u_t \ge a_1(r,t) u_{rr} + a_{21}(r,t) u_r + \frac{a_{22}(r,t)}{r} \cdot u_r - \kappa u
	\qquad \mbox{for all $r\in (0,R)$ and } t\in (0,\tm),
  \eas
  so that for all $\eps>0$, writing $\varphi(r,t):=e^{(\kappa+\eps)t} u(r,t)$ we see that
  \bea{23.4}
	\hspace*{-18mm}
	\varphi_t &\ge& e^{(\kappa+\eps)t} \cdot \Big\{
	a_1(r,t) u_{rr} + a_{21}(r,t) u_r + a_{22}(r,t) u_r - \kappa u + (\kappa+\eps)u \Big\} \nn\\
	&=& a_1(r,t) \varphi_{rr} + a_{21}(r,t) \varphi_r + a_{22}(r,t) \varphi_r +\eps e^{(\kappa+\eps)t} u
	\qquad \mbox{for all $r\in (0,R)$ and } t\in (0,\tm).
  \eea
  Now if for some $T\in (0,\tm)$, $\varphi$ attains its minimum over $[0,R] \times [0,T]$ at some
  $(r_0,t_0) \in [0,R] \times [0,T]$, then necessarily
  \be{23.5}
	\varphi_r(r_0,t_0)=0, \quad \varphi_{rr}(r_0,t_0) \ge 0 \quad \mbox{and} \quad \varphi_t(r_0,t_0) \le 0.
  \ee
  Therefore, in the case $t_0>0$ and $r_0>0$ we may directly apply (\ref{23.4}) to obtain
  \bas
	0 \ge \varphi_t(r_0,t_0)
	&\ge& a_1(r_0,t_0) \varphi_{rr}(r_0,t_0) + a_{21}(r_0,t_0) \varphi_r(r_0,t_0)
	+ \frac{a_{22}(r_0,t_0)}{r_0} \cdot \varphi_r(r_0,t_0)
	+ \eps e^{(\kappa+\eps)t_0} u(r_0,t_0) \\
	&\ge& \eps e^{(\kappa+\eps)t_0} u(r_0,t_0)
	>0,
  \eas
  which is impossible.

  However, if $t_0>0$ and $r_0=0$, then there must exist a sequence
  $(r_j)_{j\in\N}$ of numbers $r_j\in (0,R)$ such that $r_j\searrow 0$
  as $j\to\infty$ and $\varphi_r(r_j,t_0)\ge 0$ for all $j\in\N$, because otherwise $\varphi(\cdot,t_0)$
  would have a strict local maximum at $r=0$.
  Since $a_3\ge 0$, evaluating (\ref{23.4}) at $r=r_j=$ we would thus obtain that
  \bas
	\varphi_t(r_j,t_0)
	&\ge& a_1(r_j,t_0) \varphi_{rr}(r_j,t_0) + a_{21}(r_j,t_0) \varphi_r(r_j,t_0)
	+ \frac{a_{22}(r_j,t_0)}{r_j} \cdot \varphi_r(r_j,t_0) + \eps e^{(\kappa+\eps)t_0} u(r_j,t_0) \\
	&\ge& a_1(r_j,t_0) \varphi_{rr}(r_j,t_0) + a_{21}(r_j,t_0) \varphi_r(r_j,t_0)
	+ \eps e^{(\kappa+\eps)t_0} u(r_j,t_0)
	\qquad \mbox{for all } j\in\N,
  \eas
  so that since $\varphi(\cdot,t_0)$ is smooth in $[0,R]$,
  we may let $j\to\infty$ here to infer using (\ref{23.5}) that
  \bas
	0 \ge \varphi_t(0,t_0)
	&\ge& a_1(0,t_0) \varphi_{rr}(0,t_0) + a_{21}(0,t_0) \varphi_r(0,t_0)
	+ \eps e^{(\kappa+\eps)t_0} u(0,t_0) \\
	&\ge& \eps e^{(\kappa+\eps)t_0} u(0,t_0)
	>0.
  \eas
  This absurd conclusion shows that actually $t_0=0$, which implies that
  $\varphi\ge \inf_{r\in (0,R)} \varphi(r,0)=\inf_{r\in (0,R)} u_0(r)$
  throughout $[0,R] \times [0,T]$ for any $T\in (0,\tm)$.
  Taking $T\nearrow \tm$ and $\eps\searrow 0$ we thereby obtain (\ref{23.1}).
\qed

\mysection{A pointwise lower estimate for $u_r$}\label{sect4}
It remains to exclude the possibility of finite-time blow-up of $u_r$ despite boundedness of $u$.
A first step toward this can accomplished by invoking parabolic comparison to derive the following lower bound for
$u_r$ from (\ref{35.01}).
Let us emphasize that our argument makes essential use of the fact that on the right-hand side of (\ref{35.1}),
the most singular term $\frac{n-1}{r^2} \cdot \frac{uu_r}{\sqrt{u^2+u_r^2}}$ therein
appears with a negative sign, and that
in consequence a corresponding upper estimate for $u_r$ can apparently not be obtained by
a direct approach of the type pursued here, at least not when $n\ge 2$.
\begin{lem}\label{lem34}
  Assume that $\tm<\infty$, but that $\sup_{(r,t)\in (0,R)\times (0,\tm)} u(r,t)<\infty$.
  Then there exists $C>0$ such that
  \be{34.1}
	u_r(r,t) \ge -C
	\qquad \mbox{for all $r\in (0,R)$ and } t\in (0,\tm).
  \ee
\end{lem}
\proof
  According to our hypothesis, we can find $c_1>0$ such that
  \be{34.2}
	u(r,t)\le c_1
	\qquad \mbox{for all $r\in (0,R)$ and } t\in (0,\tm),
  \ee
  so that Lemma \ref{lem333} provides $c_2>0$ and $c_3>0$ such that
  \be{34.3}
	|v_r(r,t)| \le c_2 r
	\quad \mbox{and} \quad
	|v_{rr}(r,t)| \le c_3
	\qquad \mbox{for all $r\in (0,R)$ and } t\in (0,\tm).
  \ee
  We now take $D\ge 1$ and $\alpha>0$ large enough fulfilling
  \be{34.4}
	u_{0r}(r) > -D
	\qquad \mbox{for all } r\in (0,R)
  \ee
  and
  \be{34.5}
	\alpha>c_4 + \frac{c_5}{D}\va
  \ee
  where
  \be{34.54}
	c_4:=2c_1\chi + c_3\chi + c_2^2 c_3 \chi R^2 + (n-1) c_2^3 \chi R^2
  \ee
  and
  \be{34.55}
	c_5:=3c_1 c_2 c_3 \chi\mu + 3c_1^2 c_2 c_3 \chi R + (n-1) c_1 c_2^3 \chi R + 3(n-1) c_1 c_2^2 c_3 \chi R,
  \ee
  and define a comparison function $\uphi$ by letting
  \bas
	\uphi(r,t):=-D e^{\alpha t}
	\qquad \mbox{for $r\in [0,R]$ and } t\ge 0.
  \eas
  Then since $\uphi_r=\uphi_{rr}\equiv 0$, with $\parab$ as in (\ref{parab}) we have
  \bea{34.6}
	(\parab \uphi)(r,t)
	&=& -\alpha D e^{\alpha t} \nn\\
	& & - 3\frac{u_r^4 \cdot D e^{\alpha t}}{\sqrt{u^2+u_r^2}^5}
	- \frac{n-1}{r^2} \cdot \frac{u\cdot D e^{\alpha t}}{\sqrt{u^2+u_r^2}}
	- \chi\mu \frac{D e^{\alpha t}}{\sqrt{1+v_r^2}^3} \nn\\
	& & + 2\chi \frac{u\cdot D e^{\alpha t}}{\sqrt{1+v_r^2}^3}
	- \chi \frac{v_{rr} \cdot D e^{\alpha t}}{\sqrt{1+v_r^2}}
	+ \chi \frac{v_r^2 v_{rr} \cdot D e^{\alpha t}}{\sqrt{1+v_r^2}^3} \nn\\
	& & - \chi \cdot \frac{n-1}{r} \cdot \frac{v_r^3 \cdot D e^{\alpha t}}{\sqrt{1+v_r^2}^3} \nn\\
	& & - \frac{n-1}{r} \cdot \frac{u_r^4}{\sqrt{u^2+u_r^2}^3}
	- 3\chi\mu \frac{u v_r v_{rr}}{\sqrt{1+v_r^2}^5}
	+ 3\chi \frac{u^2 v_r v_{rr}}{\sqrt{1+v_r^2}^5} \nn\\
	& & - \chi \cdot \frac{n-1}{r^2} \cdot \frac{uv_r^3}{\sqrt{1+v_r^2}^3}
	+ 3\chi \cdot \frac{n-1}{r} \cdot \frac{uv_r^2 v_{rr}}{\sqrt{1+v_r^2}^5}
  \eea
  for all $r\in (0,R)$ and $t\in (0,\tm)$.
  Here the second, third, fourth and ninth term on the right are nonpositive, and we claim that each of
  the remaining summands containing $\chi$ can be controlled in modulus by the first term on the right-hand side
  suitably.
  Indeed, repeatedly using (\ref{34.2}), (\ref{34.3}) and (\ref{34.4}), we can estimate
  \bas
	\bigg| 2\chi \frac{u\cdot D e^{\alpha t}}{\sqrt{1+v_r^2}^3} \bigg|
	\le 2\chi \cdot c_1 \cdot D e^{\alpha t},
  \eas
  \bas
	\bigg| -\chi \frac{v_{rr} \cdot D e^{\alpha t}}{\sqrt{1+v_r^2}} \bigg|
	\le \chi \cdot c_3 \cdot D e^{\alpha t},
  \eas
  \bas
	\bigg| \chi \frac{v_r^2 v_{rr} \cdot D e^{\alpha t}}{\sqrt{1+v_r^2}^3} \bigg|
	\le \chi \cdot c_2^2 r^2 \cdot c_3 \cdot D e^{\alpha t}
	\le c_2^2 c_3 \chi R^2 \cdot D e^{\alpha t},
  \eas
  \bas
	\bigg| -3\chi\mu \frac{uv_r v_{rr}}{\sqrt{1+v_r^2}^5} \bigg|
	\le 3\chi\mu \cdot c_1 \cdot c_2 r \cdot c_3
	\le 3c_1c_2c_3 \chi\mu R
  \eas
    as well as
  \bas
	\bigg| 3\chi \frac{u^2 v_r v_{rr}}{\sqrt{1+v_r^2}^5} \bigg|
	\le 3\chi \cdot c_1^2 \cdot c_2 r \cdot c_3
	\le 3c_1^2 c_2 c_3 \chi R
  \eas
  and finally,
  \bas
	\bigg| -\chi \cdot \frac{n-1}{r^2} \cdot \frac{uv_r^3}{\sqrt{1+v_r^2}^3} \bigg|
	\le \chi \cdot \frac{n-1}{r^2} \cdot c_1 \cdot c_2^3 r^3
	\le (n-1) c_1 c_2^3 \chi R,
  \eas
  as well as
  \bas
	\bigg| - \chi \cdot \frac{n-1}{r} \cdot \frac{v_r^3 \cdot D e^{\alpha t}}{\sqrt{1+v_r^2}^3} \bigg|
	\le \chi \cdot \frac{n-1}{r} \cdot c_2^3 r^3 \cdot D e^{\alpha t}
	\le (n-1) c_2^3 \chi R^2 \cdot D e^{\alpha t}
  \eas
  and
  \bas
	\bigg| 3\chi \cdot \frac{n-1}{r} \cdot \frac{uv_r^2 v_{rr}}{\sqrt{1+v_r^2}^5} \bigg|
	\le 3\chi \cdot \frac{n-1}{r} \cdot c_1 \cdot c_2^2 r^2 \cdot c_3
	\le 3(n-1) c_1 c_2^2 c_3 \chi R.
  \eas

  Therefore, (\ref{34.6}) implies that with $c_4$ and $c_5$ as in (\ref{34.54}) and (\ref{34.55}) we have
  \bas
	(\parab \uphi)(r,t)
	&\le& -\alpha D e^{\alpha t} + c_4 \cdot D e^{\alpha t} + c_5 \\
	&\le& -\alpha D e^{\alpha t} + \Big(c_4+\frac{c_5}{D}\Big) \cdot D e^{\alpha t}
	\qquad \mbox{for all $r\in (0,R)$ and } t\in (0,\tm),
  \eas
  whence our assumption (\ref{34.5}) on $\alpha$ ensures that $(\parab \uphi)(r,t)< 0$
  for all $r\in (0,R)$ and $t\in (0,\tm)$.
  Since $(\parab u_r)(r,t)=0$ for all $(r,t)\in (0,R)\times (0,\tm)$ by Lemma \ref{lem35}, and since moreover
  \bas
	\uphi(r,0)=-D<u_{0r}(r)=u_r(0,r)
	\qquad \mbox{for all } r\in [0,R]
  \eas
  and, clearly,
  \bas
	\uphi_r(0,t)=u_r(0,t)=0
	\quad \mbox{as well as} \quad
	\uphi_r(R,t)=u_r(R,t)=0
	\qquad \mbox{for all } t\in (0,\tm),
  \eas
  from the comparison principle we conclude that $u_r(r,t)\ge \uphi(r,t)$ for all $r\in (0,R)$ and
  $t\in (0,\tm)$, and that hence
  \bas
	u_r(r,t) \ge -D e^{\alpha \tm}
	\qquad \mbox{for all $r\in (0,R)$ and } t\in (0,\tm),
  \eas
  which proves the claim.
\qed
\mysection{A bound for $|u_r|$. Proof of Theorem \ref{theo43}}\label{sect5}	
%
%
%
%
%
The goal of this section is to complete the proof of Theorem \ref{theo43} by
further developing the one-sided inequality for $u_r$ from Lemma \ref{lem34}
into a bound for $|u_r|$ in modulus, provided that $\tm$ is finite but $u$ itself remains bounded
(Corollary \ref{cor42}).
An important role in our analysis in this direction	
will be played by the function $z:=\frac{u_t}{u}$,
which is indeed well-defined and continuous in $[0,R]\times [0,\tm)$ by Lemma \ref{lem21}.
Furthermore, according to Lemma \ref{lem36} we have the representation
\bea{z}
	z &=& \frac{u^2 u_{rr}}{\sqrt{u^2+u_r^2}^3} + \frac{u_r^4}{u\sqrt{u^2+u_r^2}^3}
		+ \frac{n-1}{r} \cdot \frac{u_r}{\sqrt{u^2+u_r^2}} \nn\\
	& & - \chi \frac{u_r v_r}{u\sqrt{1+v_r^2}} - \chi \frac{\mu-u}{\sqrt{1+v_r^2}^3}
	- \chi \cdot \frac{n-1}{r} \cdot \frac{v_r^3}{\sqrt{1+v_r^2}^3}\va
\eea
for $r\in (0,R)$ and $t\in (0,\tm)$.\abs
Now in a first key observation, to be presented in Corollary \ref{cor40},
we will establish a useful relationship between $u_r$ and the function $z$,
essentially controlling $\|u_r(\cdot,t)\|_{L^\infty((0,R))}$ for any fixed $t\in (0,\tm)$
by the maximum of the {\em positive part} $z_+$ of $z$ over the {\em whole memory region} $(0,R)\times (0,t)$.
To achieve this, we will foremost use an integral technique to estimate $|u_r|$ in terms of $z_+$ on the basis
of (\ref{z}) and Lemma \ref{lem34} in a suitably small subinterval $(0,R_0)$ of $(0,R)$ (Lemma \ref{lem32}).
This will in particular imply an upper bound for $u_r$ at $r=R_0$ and therefore allow for applying
a comparison argument to (\ref{35.02})
which will yield a pointwise upper estimate for $u_r$ in the corresponding outer region $(R_0,R)$
(Lemma \ref{lem37}).\abs
The second essential step will thereafter consist in deriving a nonlocal parabolic inequality for $z$
with a memory-type nonlinearity (Lemma \ref{lem39}).
Upon another comparison, this will entail a pointwise upper bound for $z$ (Lemma \ref{lem41})
and hence also for $u_r$.
\subsection{A bound for $|u_r|$ in terms of $z_+$}
\subsubsection{Estimating $|u_r|$ near the origin}
Let us first apply an appropriate testing procedure to (\ref{z}) to find some small $R_0\in (0,R)$ with the
property that $u_r(\cdot,t)$ can be bounded in certain weighted Lebesgue spaces over $(0,R_0)$
in such a way that on taking limits we can derive a respective $L^\infty$ estimate from this.

\begin{lem}\label{lem32}
  Assume that $\tm<\infty$, but that $\sup_{(r,t)\in (0,R)\times (0,\tm)} u(r,t)<\infty$.
  Then there exist $R_0\in (0,R)$ and $C>0$ such that
  \be{32.01}
	\|u_r(\cdot,t)\|_{L^\infty((0,R_0))} \le C \cdot \Big(1+\|z_+(\cdot,t)\|_{L^\infty((0,R_0))} \Big)
	\qquad \mbox{for all } t\in (0,\tm).
  \ee
\end{lem}
\proof
  We first rearrange (\ref{z}) to obtain
  \bea{32.1}
	\frac{u_r^4}{u^3} &=& \frac{\sqrt{u^2+u_r^2}}{u^2} \cdot z
	- u_{rr} - \frac{n-1}{r} \cdot \frac{\sqrt{u^2+u_r^2}^2 u_r}{u^2} \nn\\
	& & + \chi \frac{(\mu-u)\sqrt{u^2+u_r^2}^3}{\sqrt{1+v_r^2}^3}
	+ \chi \frac{\sqrt{u^2+u_r^2}^3 u_r v_r}{u^3 \sqrt{1+v_r^2}}
	+ \chi \cdot \frac{n-1}{r} \cdot \frac{\sqrt{u^2+u_r^2}^3 v_r^3}{u^2 \sqrt{1+v_r^2}^3},
  \eea
  where
  \bea{32.100}
	- u_{rr} - \frac{n-1}{r} \cdot \frac{\sqrt{u^2+u_r^2}^2 u_r}{u^2}
	&=& - \Big(u_{rr}+\frac{n-1}{r} u_r\Big) - \frac{n-1}{r} \cdot \frac{u_r^3}{u^2} \nn\\
	&=& - \frac{1}{r^{n-1}} (r^{n-1} u_r)_r
	-\frac{n-1}{r} \cdot \frac{u_r^3}{u^2}
  \eea
  for $r\in (0,R)$ and $t\in (0,\tm)$.
  In order to choose $R_0$ appropriately, we use our boundedness assumption on $u$ to fix $c_1\ge \mu$ and
  $c_2>0$ such that
  \be{32.2}
	u(r,t) \le c_1
	\qquad \mbox{for all $r\in (0,R)$ and } t\in (0,\tm)
  \ee
  and
  \be{32.3}
	\sqrt{u^2+u_r^2}^3 \le c_2 \cdot (1+|u_r|^3)
	\qquad \mbox{in } (0,R)\times (0,\tm),
  \ee
  and recall Lemma \ref{lem23} to find $c_3>0$ fulfilling
  \be{32.4}
	u(r,t)\ge c_3
	\qquad \mbox{for all $r\in (0,R)$ and } t\in (0,\tm).
  \ee
  We claim that then the conclusion of the lemma holds if we pick any $R_0\in (0,R)$ satisfying
  \be{32.5}
	R_0 \le \frac{n c_3^3}{4c_1^3 c_2 \chi\mu}.
  \ee
  To see this, we take an arbitrary even integer $m\ge 0$, multiply (\ref{32.1}) by $r^{n-1} u_r^m$
  and integrate over $(0,R_0)$ to see using (\ref{32.100}) that
  \bea{32.6}
	I(t) &:=&
	\int_0^{R_0} r^{n-1} \frac{u_r^{m+4}}{u^3} dr \nn\\
	&=& \int_0^{R_0} r^{n-1} \frac{\sqrt{u^2+u_r^2}^3}{u^2} \cdot u_r^m z dr
	- \int_0^{R_0} (r^{n-1} u_r)_r \cdot u_r^m dr
	- (n-1) \int_0^{R_0} r^{n-2} \frac{u_r^{m+3}}{u^2} dr \nn\\
	& & + \chi \int_0^{R_0} r^{n-1} \frac{(\mu-u) \sqrt{u^2+u_r^2}^3 u_r^m}{u^2 \sqrt{1+v_r^2}^3} dr
	+ \chi \int_0^{R_0} r^{n-1} \frac{\sqrt{u^2+u_r^2}^3 u_r^{m+1} v_r}{u^3\sqrt{1+v_r^2}} dr \nn\\
	& & + (n-1)\chi \int_0^{R_0} r^{n-2} \frac{\sqrt{u^2+u_r^2}^3 u_r^m v_r^3}{u^2 \sqrt{1+v_r^2}^3} dr
	\nn\\[2mm]
	&=:& J_1(t)+...+J_6(t)
	\qquad \mbox{for all } t\in (0,\tm).
  \eea
  Here by (\ref{32.2}) we have
  \be{32.7}
	I(t) \ge \frac{1}{c_1^3} \cdot \int_0^{R_0} r^{n-1} u_r^{m+4} dr
	\qquad \mbox{for all } t\in (0,\tm),
  \ee
  and our goal is to show that the sum on the right-hand side of (\ref{32.6}) can be controlled adequately by
  the term on the right of (\ref{32.7}).\\
  For this purpose, we first use Lemma \ref{lem33} in rewriting $J_5(t)$ to obtain
  \bea{32.8}
	J_5(t)
	&=& \frac{\chi \mu}{n} \int_0^{R_0} r^n \cdot \frac{\sqrt{u^2+u_r^2}^3 u_r^{m+1}}{u^3 \sqrt{1+v_r^2}} dr
	- \chi \int_0^{R_0} \frac{\sqrt{u^2+u_r^2}^3 u_r^{m+1}}{u^3 \sqrt{1+v_r^2}}
		\cdot \bigg( \int_0^r \rho^{n-1} u(\rho,t) d\rho \bigg) dr \nn\\
	&=:& J_{51}(t)+J_{52}(t)
	\qquad \mbox{for all } t\in (0,\tm),
  \eea
  where (\ref{32.3}), (\ref{32.4}) and Young's inequality enable us to infer that
  \bas
	J_{51}(t)
	&\le& \frac{\chi\mu}{n} \int_0^{R_0} r^n \cdot \frac{c_2 \cdot (1+|u_r|^3) \cdot |u_r|^{m+1}}{c_3^3} dr \\
	&=& \frac{c_2 \chi\mu}{nc_3^3} \int_0^{R_0} r^n |u_r|^{m+1} dr
	+ \frac{c_2 \chi\mu}{nc_3^3} \int_0^{R_0} r^n u_r^{m+4} dr \\
	&\le& \frac{c_2\chi\mu}{nc_3^3} \int_0^{R_0} r^n (1+u_r^{m+4}) dr
	+ \frac{c_2 \chi\mu}{nc_3^3} \int_0^{R_0} r^n u_r^{m+4} dr \\
	&=& \frac{c_2 \chi\mu R_0^{n+1}}{n(n+1)c_3^3}
	+ \frac{2c_2 \chi\mu}{nc_3^3} \int_0^{R_0} r^n u_r^{m+4}dr
	\qquad \mbox{for all  } t\in (0,\tm).
  \eas
  Trivially estimating
  \bas
	\int_0^{R_0} r^n u_r^{m+4} dr
	\le R_0 \int_0^{R_0} r^{n-1} u_r^{m+4},
  \eas
  according to (\ref{32.7}) and our restriction (\ref{32.5}) on $R_0$ we thus conclude that
  \be{32.9}
	J_{51}(t)
	\le c_4 + \frac{1}{2c_1^3} \int_0^{R_0} r^{n-1} u_r^{m+4} dr
	\qquad \mbox{for all } t\in (0,\tm)
  \ee
  with $c_4:=\frac{c_2 \chi\mu R^{n+1}}{n(n+1)c_3^3}$. \\
  In order to derive an appropriate upper bound for the second term on the right of (\ref{32.8}), let us apply
  Lemma \ref{lem34} to fix a constant $L\ge 1$ such that
  \be{32.10}
	u_r(r,t) \ge -L
	\qquad \mbox{for all $r\in (0,R)$ and } t\in (0,\tm).
  \ee
  Then since $-u_r^{m+1} = (-u_r)^{m+1}$ due to the fact that $m+1$ is odd, this in conjunction with (\ref{32.3}),
  (\ref{32.4}) and (\ref{32.1}) implies that
  \bea{32.11}
	J_{52}(t)
	&\le& \chi L^{m+1} \int_0^{R_0} \frac{\sqrt{u^2+u_r^2}^3}{u^3 \sqrt{1+v_r^2}}
	\cdot \bigg(\int_0^r \rho^{n-1} u(\rho,t) d\rho\bigg) dr \nn\\
	&\le& \frac{c_1 c_2 \chi}{c_3^3} L^{m+1} \int_0^{R_0} (1+|u_r|^3)
	\cdot \Big( \int_0^r \rho^{n-1} d\rho\Big) dr \nn\\
	&=& \frac{c_1 c_2 \chi}{nc_3^3} L^{m+1} \int_0^{R_0} r^n (1+|u_r|^3) dr \nn\\
	&=& \frac{c_1 c_2 \chi R_0^{n+1}}{n(n+1) c_3^3} L^{m+1}
	+ \frac{c_1 c_2 \chi}{nc_3^3} \int_0^{R_0} r^n L^{m+1} |u_r|^3 dr
	\qquad \mbox{for all } t\in (0,\tm).
  \eea
  As by Young's inequality we have
  \bas
	\int_0^{R_0} r^n L^{m+1} |u_r|^3 dr
	&\le& \frac{3}{m+4} \int_0^{R_0} r^n u_r^{m+4} dr
	+ \frac{m+1}{m+4} \int_0^{R_0} r^n L^{m+4} dr \\
	&\le& \frac{3R}{m+4} \int_0^{R_0} r^{n-1} u_r^{m+4} dr
	+ \frac{R^{n+1}}{n+1} L^{m+4},
  \eas
  in view of the fact that $L\ge 1$ we obtain from (\ref{32.11}) that
  \be{32.12}
	J_{52}(t)
	\le c_5 L^{m+4} + \frac{c_6}{m+4} \int_0^{R_0} r^{n-1} u_r^{m+4} dr
	\qquad \mbox{for all } t\in (0,\tm)
  \ee
  with $c_5:=\frac{2c_1 c_2 \chi R^{n+1}}{n(n+1)c_3^3}$ and $c_6:=\frac{3c_1 c_2 \chi R}{nc_3^3}$. \\
  Going back to (\ref{32.6}), we next apply (\ref{32.3}), (\ref{32.4}) and Young's and H\"older's inequalities
  to estimate
  \bas
	J_1(t)
	&\le& \frac{c_2}{c_3^2} \int_0^{R_0} r^{n-1} (1+|u_r|^3) \cdot u_r^m z_+ dr \\
	&=& \frac{c_2}{c_3^2} \int_0^{R_0} r^{n-1} u_r^m z_+ dr
	+ \frac{c_2}{c_3^2} \int_0^{R_0} r^{n-1} |u_r|^{m+3} z_+ dr \\
	&\le& \frac{c_2}{c_3^2} \int_0^{R_0} r^{n-1} z_+ dr
	+ \frac{2c_2}{c_3^2} \int_0^{R_0} r^{n-1} |u_r|^{m+3} z_+ dr \\
	&\le& \frac{c_2}{c_3^2} \bigg( \int_0^{R_0} r^{n-1} z_+^{m+4} dr \bigg)^\frac{1}{m+4}
		\cdot \Big(\frac{R_0^n}{n}\Big)^\frac{m+3}{m+4}
	+\frac{2c_2}{c_3^2} \bigg(\int_0^{R_0} r^{n-1} u_r^{m+4} dr \bigg)^\frac{m+3}{m+4}
		\cdot \bigg(\int_0^{R_0} r^{n-1} z_+^{m+4} dr\bigg)^\frac{1}{m+4},
  \eas
  which means that if we let $R_1:=\max \{1,R\}$ and $c_7:=\max \{ \frac{c_2}{c_3^2} \cdot \frac{R_1^n}{n} \, , \,
  \frac{2c_2}{c_3^2} \}$, then
  \be{32.13}
	J_1(t) \le c_7 \cdot \Bigg\{ 1 + \bigg( \int_0^{R_0} r^{n-1} u_r^{m+4}dr \bigg)^\frac{m+3}{m+4} \Bigg\}
	\cdot \bigg( \int_0^{R_0} r^{n-1} z_+^{m+4} dr \bigg)^\frac{1}{m+4}
	\qquad \mbox{for all } t\in (0,\tm).
  \ee
  As for the second term on the right of (\ref{32.6}), in order to remove second-order derivatives, in the case
  $n\ge 2$ we twice integrate by parts to obtain
  \bas
	J_2(t)
	&=& m \int_0^{R_0} r^{n-1} u_r^m u_{rr} dr
	- R_0^{n-1} u_r^{m+1} (R_0,t) \\
	&=& - \frac{(n-1)m}{m+1} \int_0^{R_0} r^{n-2} u_r^{m+1} dr
	+ \frac{m}{m+1} R_0^{n-1} u_r^{m+1}(R_0,t)
	- R_0^{n-1} u_r^{m+1} (R_0,t) \\
	&=& - \frac{(n-1)m}{m+1} \int_0^{R_0} r^{n-2} u_r^{m+1} dr
	- \frac{1}{m+1} R_0^{n-1} u_r^{m+1}(R_0,t)
	\qquad \mbox{for all } t\in (0,\tm).
  \eas
  Once more since $m+1$ is odd, (\ref{32.10}) again becomes applicable to provide the one-sided estimate
  \bea{32.14}
	J_2(t)
	&\le& \frac{(n-1)m}{m+1} L^{m+1} \int_0^{R_0} r^{n-2}
	+\frac{1}{m+1} R_0^{n-1} L^{m+1} \nn\\
	&\le& c_8 L^{m+1}
	\qquad \mbox{for all } t\in (0,\tm)
  \eea
  with $c_8:=\frac{(n-1)m}{m+1} \cdot \frac{R^{n-1}}{n-1} + \frac{1}{m+1} R^{n-1} \equiv R^{n-1}$
  when $n\ge 2$, and it can easily be verified that this conclusion can be extended so as to include the case $n=1$
  as well.\\

  Similarly, since also $m+3$ is odd, we may invoke (\ref{32.10}) and then (\ref{32.4}) to see that in the case
  $n\ge 2$,
  \bea{32.15}
	J_3(t)
	&\le& (n-1)\int_0^{R_0} r^{n-2} \frac{L^{m+3}}{u^2} dr \nn\\
	&\le& \frac{n-1}{c_3^2} L^{m+3} \frac{R_0^{n-1}}{n-1} \nn\\
	&\le& c_9 L^{m+3}
	\qquad \mbox{for all } t\in (0,\tm)
  \eea
  with $c_9:=\frac{R^{n-1}}{c_3^2}$, and note that (\ref{32.15}) trivially holds when $n=1$.\\

  We next estimate $J_4(t)$ by first using (\ref{32.3}), (\ref{32.4}) and Young's inequality according to
  \bas
	J_4(t)
	&\le& \frac{c_2 \chi\mu}{c_3^2} \int_0^{R_0} r^{n-1} (1+|u_r|^3) u_r^m dr \\
	&\le& \frac{c_2\chi\mu}{c_3^2} \int_0^{R_0} r^{n-1} dr
	+ \frac{2c_2\chi\mu}{c_3^2} \int_0^{R_0} r^{n-1} |u_r|^{m+3} dr \\
	&=& \frac{c_2\chi\mu R_0^n}{nc_3^2}
	+ \frac{2c_2\chi\mu}{c_3^2} \int_0^{R_0} r^{n-1} |u_r|^{m+3} dr
	\qquad \mbox{for all } t\in (0,\tm).
  \eas
  Since thanks to the H\"older inequality we know that
  \bea{32.155}
	\int_0^{R_0} r^{n-1} |u_r|^{m+3} dr
	&\le& \Big(\frac{R_0^n}{n}\Big)^\frac{1}{m+4} \cdot
		\bigg( \int_0^{R_0} r^{n-1} u_r^{m+4} dr \bigg)^\frac{m+3}{m+4} \nn\\
	&\le& \frac{R_1^n}{n} \cdot
		\bigg( \int_0^{R_0} r^{n-1} u_r^{m+4} dr \bigg)^\frac{m+3}{m+4},
  \eea
  again with $R_1=\max \{1,R\}$, this entails that
  \be{32.16}
	J_4(t) \le c_{10} \cdot \Bigg\{
	1 + \bigg( \int_0^{R_0} r^{n-1} u_r^{m+4} dr \bigg)^\frac{m+3}{m+4} \Bigg\}
	\qquad \mbox{for all } t\in (0,\tm)
  \ee
  if we let
  $c_{10}:=\max\{ \frac{c_2\chi\mu R^n}{nc_3^2} \, , \, \frac{2c_2\chi\mu}{c_3^2} \cdot \frac{R_1^n}{n} \}$.\\

  Finally, in treating the last integral in (\ref{32.6}) we make use of the upper estimate for $v_r$ in
  (\ref{333.0}) to see, again recalling (\ref{32.3}), (\ref{32.4}) and using Young's inequality, that
  \bas
	J_6(t)
	&\le& \frac{(n-1)\chi\mu^3}{n^3} \int_0^{R_0} \frac{r^{n+1} \sqrt{u^2+u_r^2}^3 u_r^m}{u^2\sqrt{1+v_r^2}^3}
		dr \\
	&\le& \frac{(n-1)c_2 \chi \mu^3}{n^3 c_3^2} \int_0^{R_0} r^{n+1} (1+|u_r|^3) u_r^m dr \\
	&\le& \frac{(n-1)c_2 \chi \mu^3}{n^3 c_3^2} \int_0^{R_0} r^{n+1} dr
	+ \frac{2(n-1)c_2 \chi \mu^3}{n^3 c_3^2} \int_0^{R_0} r^{n+1} |u_r|^{m+3} dr \\
	&=& \frac{(n-1)c_2 \chi \mu^3 R_0^{n+2}}{n^3 (n+2) c_3^2}
	+ \frac{2(n-1)c_2 \chi \mu^3}{n^3 c_3^2} \int_0^{R_0} r^{n+1} |u_r|^{m+3} dr
	\qquad \mbox{for all } t\in (0,\tm).
  \eas
  Now due to (\ref{32.155}),
  \bas
	\int_0^{R_0} r^{n+1} |u_r|^{m+3} dr
	&\le& R_0^2 \int_0^{R_0} r^{n-1} |u_r|^{m+3} dr \\
	&\le& \frac{R_1^{n+2}}{n} \bigg(\int_0^{R_0} r^{n-1} u_r^{m+4} dr \bigg)^\frac{m+3}{m+4},
  \eas
  whence we obtain that
  \be{32.17}
	J_6(t) \le c_{11} \cdot \Bigg\{ 1 + \bigg( \int_0^{R_0} r^{n-1} u_r^{m+4} dr \bigg)^\frac{m+3}{m+4} \Bigg\}
	\qquad \mbox{for all } t\in (0,\tm)
  \ee
  with $c_{11}:=\max\Big\{ \frac{(n-1) c_2 \chi \mu^3 R^{n+2}}{n^3(n+2)c_3^2} \, , \,
  \frac{2(n-1) c_2 \chi\mu^3}{n^3 c_2^3} \cdot \frac{R_1^{n+2}}{n} \Big\}$.\abs
  In summary, (\ref{32.7}), (\ref{32.13}), (\ref{32.14}), (\ref{32.15}), (\ref{32.8}), (\ref{32.9}), (\ref{32.12})
  and (\ref{32.17}) combined with (\ref{32.6}) show that
  \bas
	\frac{1}{c_1^3} \int_0^{R_0} r^{n-1} u_r^{m+4} dr
	&\le& c_7 \cdot \Bigg\{ 1+\bigg(\int_0^{R_0} r^{n-1} u_r^{m+4} dr \bigg)^\frac{m+3}{m+4} \Bigg\}
		\cdot \bigg(\int_0^{R_0} r^{n-1} z_+^{m+4} dr \bigg)^\frac{1}{m+4} \\
	& & + c_8 L^{m+1} + c_9 L^{m+3} \\
	& & + c_{10} \cdot \Bigg\{ 1+\bigg(\int_0^{R_0} r^{n-1} u_r^{m+4} dr \bigg)^\frac{m+3}{m+4} \Bigg\} \\
	& & + c_4 + \frac{1}{2c_1^3} \int_0^{R_0} r^{n-1} u_r^{m+4} dr \\
	& & + c_5 L^{m+4} + \frac{c_6}{m+4} \int_0^{R_0} r^{n-1} u_r^{m+4} dr \\
	& & + c_{11} \Bigg\{ 1+\bigg(\int_0^{R_0} r^{n-1} u_r^{m+4} dr \bigg)^\frac{m+3}{m+4} \Bigg\}
	\qquad \mbox{for all } t\in (0,\tm).
  \eas
  Since $L\ge 1$, this means that if $m_\star$ is sufficiently large such that
  \be{32.18}
	\frac{c_6}{m_\star+4} \le \frac{1}{4c_1^3},
  \ee
  and if $m\ge m_\star$, then
  \bea{32.19}
	\frac{1}{4c_1^3} \int_0^{R_0} r^{n-1} u_r^{m+4} dr
	&\le& c_7 \cdot \Bigg\{ 1+\bigg(\int_0^{R_0} r^{n-1} u_r^{m+4} dr \bigg)^\frac{m+3}{m+4} \Bigg\}
	\cdot \bigg( \int_0^{R_0} r^{n-1} z_+^{m+4} dr \bigg)^\frac{1}{m+4} \nn\\
	& & + c_{12} \bigg(\int_0^{R_0} r^{n-1} u_r^{m+4} dr\bigg)^\frac{m+3}{m+4} \nn\\
	& & + c_{13} L^{m+4}
	\qquad \mbox{for all } t\in (0,\tm),
  \eea
  where $c_{12}:=c_{10}+c_{11}$ and $c_{13}:=c_8+c_9+c_{10}+c_4+c_5+c_{11}$.\abs
  Now in order to prove (\ref{32.01}) for some suitably large $C>0$ independent of $t\in (0,\tm)$, we fix
  any such $t$ and first consider the case when there exists a sequence of even numbers $m=m_j\ge m_\star, j\in\N$,
  such that $m_j\to\infty$ as $j\to\infty$ and
  \be{32.20}
	\bigg(\int_0^{R_0} r^{n-1} u_r^{m+4} dr \bigg)^\frac{m+3}{m+4}
	\le L^{m+4}
	\qquad \mbox{for all } m\in (m_j)_{j\in\N}.
  \ee
  Then taking $j\to\infty$ here, we directly obtain that
  \be{32.21}
	\|u_r(\cdot,t)\|_{L^\infty((0,R_0))}
	= \lim_{j\to\infty} \bigg(\int_0^{R_0} r^{n-1} u_r^{m_j+4} (r,t) dr \bigg)^\frac{1}{m_j+4}
	= \lim_{j\to\infty} L^\frac{m_j+4}{m_j+3} = L.
  \ee
  If conversely, such a sequence does not exist, then we can find $m_{\star\star} \ge m_\star$ such that for all
  even $m\ge m_{\star\star}$,
  \bas
	\bigg(\int_0^{R_0} r^{n-1} u_r^{m+4} dr \bigg)^\frac{m+3}{m+4} > L^{m+4}.
  \eas
  Using that $L\ge 1$, we thus infer from (\ref{32.19}) that for any such $m$ we have
  \bas
	\frac{1}{4c_1^3} \int_0^{R_0} r^{n-1} u_r^{m+4}(r,t) dr
	&<& 2c_7 \bigg(\int_0^{R_0} r^{n-1} u_r^{m+4}(r,t) dr\bigg)^\frac{m+3}{m+4} \cdot
	\bigg(\int_0^{R_0} r^{n-1} z_+^{m+4} (r,t)dr \bigg)^\frac{1}{m+4} \\
	& & + c_{12} \bigg(\int_0^{R_0} r^{n-1} u_r^{m+4}(r,t) dr\bigg)^\frac{m+3}{m+4}
	+ c_{13} \bigg(\int_0^{R_0} r^{n-1} u_r^{m+4}(r,t) dr\bigg)^\frac{m+3}{m+4}
  \eas
  and hence
  \bas
	\frac{1}{4c_1^3} \bigg(\int_0^{R_0} r^{n-1} u_r^{m+4}(r,t) dr\bigg)^\frac{1}{m+4}
	\le 2c_7 \bigg(\int_0^{R_0} r^{n-1} z_+^{m+4}(r,t) dr \bigg)^\frac{1}{m+4} + c_{12}+c_{13}.
  \eas
  In the limit $m\to\infty$, we therefore conclude that in this case,
  \be{32.22}
	\frac{1}{4c_1^3} \cdot \|u_r(\cdot,t)\|_{L^\infty((0,R_0))}
	\le 2c_7 \cdot \|z_+(\cdot,t)\|_{L^\infty((0,R_0))} + c_{12}+c_{13}.
  \ee
  Since $c_7, c_{12}$ and $c_{13}$ as well as $L$ are independent of $t\in (0,\tm)$, (\ref{32.21}) and (\ref{32.22})
  establish (\ref{32.01}).
\qed

\subsubsection{Estimating $|u_r|$ near the boundary}
For fixed $t\in (0,\tm)$, the above lemma in particular implies an upper bound for $u_r$ in terms of
$\|z_+\|_{L^\infty((0,t)\times (0,R_0))}$ on the lateral boundary line
$r=R_0$ of the parabolic cylinder $(R_0,R)\times (0,t)$.
This will enable us to apply a comparison argument to derive an estimate from above for $u_r$
in this region on the basis of (\ref{35.02}) to achieve the following.
\begin{lem}\label{lem37}
  Assume that $\tm<\infty$, but that $\sup_{(r,t)\in (0,R)\times (0,\tm)} u(r,t)<\infty$.
  Then with $R_0\in (0,R)$ taken from Lemma \ref{lem32}, we can find $C>0$ such that
  \be{37.1}
	\|u_r(\cdot,t)\|_{L^\infty((R_0,R))} \le C \cdot \Big(1+\|z_+\|_{L^\infty((0,R_0)\times (0,t))} \Big)
	\qquad \mbox{for all } t\in (0,\tm).
  \ee	
\end{lem}
\proof
  According to Lemma \ref{lem32}, we can pick $c_1>0$ such that
  \bas
	u_r(R_0,t) \le c_1 \cdot \Big(1+\|z_+(\cdot,t)\|_{L^\infty((0,R_0))} \Big)
	\qquad \mbox{for all } t\in (0,\tm),
  \eas
  which in particular implies that given any $t_0\in (0,\tm)$ we have
  \be{37.2}
	u_r(R_0,t) \le D_1(t_0):= c_1 \cdot \Big( 1+\|z_+\|_{L^\infty((0,R_0)\times (0,t_0))} \Big)
	\qquad \mbox{for all } t\in (0,t_0).
  \ee
  Let us next use our hypothesis and recall Lemma \ref{lem23} to pick $c_2>0$ and $c_3>0$ fulfilling
  \be{37.3}
	c_2 \le u(r,t) \le c_3
	\qquad \mbox{for all $r\in (0,R)$ and } t\in (0,\tm),
  \ee
  and apply Lemma \ref{lem333} to find $c_4>0$ and $c_5>0$ such that
  \be{37.4}
	|v_r(r,t)| \le c_4 r
	\quad \mbox{and} \quad
	|v_{rr}(r,t)| \le c_5
	\qquad \mbox{for all $r\in (0,R)$ and } t\in (0,\tm).
  \ee
  Therefore, the coefficient functions $\tilde A_3$ and $\tilde A_4$ in (\ref{qarab1}) can be estimated according to
  \be{37.5}
	\tilde A_3(r,t)
	\le c_6:=\frac{n-1}{R_0} + \frac{2\chi}{c_2^2} + \chi c_5 + \chi \cdot c_4^2 R^2 \cdot c_5
	+ (n-1)\chi \cdot c_4^3 R^2
  \ee
  and
  \bea{37.6}
	\tilde A_4(r,t)
	\le C_7 &:=& 3c_3 + \frac{n-1}{R_0^2} \cdot c_2
	+ 3\chi\mu \cdot c_3 \cdot c_4 R \cdot c_5
	+ 3\chi \cdot c_3^2 \cdot c_4 R \cdot c_5 \\
	& & + (n-1)\chi \cdot c_3 \cdot c_4^3 R
	+ 3(n-1) \cdot \cdot c_3 \cdot c_4^2 R \cdot c_5
  \eea
  for all $r\in (0,R)$ and $t\in (0,\tm)$.
  We now fix $\alpha>0$ large such that
  \be{37.7}
	\alpha>c_6+\frac{c_7}{D}
  \ee
  and, given $t_0\in (0,\tm)$, define
  \bas
	\ophi(r,t):=D \, e^{\alpha t}
	\qquad \mbox{for $r\in [R_0,R]$ and } t\in [0,t_0],
  \eas
  where
  \be{37.8}
	D:=\max \Big\{ D_1(t_0) \, , \sup_{r\in (R_0,R)} u_{0r}(r) \Big\} \, +1.
  \ee
  Then (\ref{37.2}) asserts that
  \bas
	u_r(R_0,t)\le \ophi(R_0,t)
	\qquad \mbox{for all } t\in (0,t_0),
  \eas
  whereas clearly
  \bas
	u_r(R,t)=\ophi(R,t)=0
	\qquad \mbox{for all } t\in (0,t_0)
  \eas
  and
  \bas
	u_r(r,0)=u_{0r}(r) < D =\ophi(r,0)
	\qquad \mbox{for all } r\in [R_0,R].
  \eas
  Moreover, since $\ophi$ is positive and $\ophi_r=\ophi_{rr}\equiv 0$, we may use (\ref{37.5}), (\ref{37.6}) and (\ref{37.7})
  to see that with $\qarab$ as in (\ref{qarab}) we have
  \bas
	\qarab \ophi
	&=& \ophi_t - \tilde A_3(r,t) \ophi - \tilde A_4(r,t) \\
	&=& \alpha D e^{\alpha t} - \tilde A_3(r,t) \cdot D e^{\alpha t} - \tilde A_4(r,t) \\
	&\ge& (\alpha - c_6) \cdot D e^{\alpha t} - c_7 \\
	&\ge& (\alpha - c_6) D - c_7 \\
	&>& 0
	\qquad \mbox{for all $r\in (R_0,R)$ and } t\in (0,t_0).
  \eas
  As $\qarab u_r \equiv 0$ due to Lemma \ref{lem35}, by comparison we conclude that $u_r \le \ophi$ in
  $(R_0,R)\times (0,t_0)$, which in view of (\ref{37.8}) and (\ref{37.2}) readily entails (\ref{37.1}).
\qed
\subsubsection{A bound for $|u_r|$ in the entire domain}
Let us summarize the outcome of Lemma \ref{lem32} and Lemma \ref{lem37}:
\begin{cor}\label{cor40}
  If $\tm<\infty$ but $\sup_{(r,t)\in (0,R)\times (0,\tm)} u(r,t)<\infty$,
  then there exists $C>0$ such that
  \bas
	\|u_r(\cdot,t)\|_{L^\infty((0,R))} \le C \cdot \Big(1+\|z_+\|_{L^\infty((0,R)\times (0,t))} \Big)
	\qquad \mbox{for all } t\in (0,\tm).
  \eas
\end{cor}
\proof
  We only need to combine Lemma \ref{lem32} with Lemma \ref{lem37}.
\qed
\subsection{A nonlocal parabolic inequality for $z$}
In order to bound $z$ from above, let us first identify a linear inhomogeneous parabolic equation satisfied
by this function.
\begin{lem}\label{lem38}
  The function $z=\frac{u_t}{u}$ satisfies
  \be{38.1}
	z_t=B_1(r,t)z_{rr} + B_{21}(r,t) z_r + \frac{B_{22}(r,t)}{r} z_r + B_3(r,t)z + B_4(r,t)
	\qquad \mbox{for all $r\in (0,R)$ and } t\in (0,\tm),
  \ee
  where
  \be{38.2}
	\left\{ \begin{array}{l}
  \displaystyle	B_1(r,t):=\frac{u^3}{\sqrt{u^2+u_r^2}^3}, \\[2mm]
	\displaystyle B_{21}(r,t):=2\frac{u^2 u_r}{\sqrt{u^2+u_r^2}^3} - 3\frac{u^3 u_r u_{rr}}{\sqrt{u^2+u_r^2}^5}
	+ 4\frac{u_r^3}{\sqrt{u^2+u_r^2}^3} - 3\frac{u_r^5}{\sqrt{u^2+u_r^2}^5} -\chi \frac{v_r}{\sqrt{1+v_r^2}}, \\[2mm]
\displaystyle	B_{22}(r,t):=(n-1)\frac{u^3}{\sqrt{u^2+u_r^2}^3}, \\[2mm]
\displaystyle	B_3(r,t):=\chi \frac{u}{\sqrt{1+v_r^2}^3} \qquad \qquad \mbox{and} \\[2mm]
\displaystyle	B_4(r,t):=-3\chi \frac{u(\mu-u) u_r v_r}{\sqrt{u^2+u_r^2} \cdot \sqrt{1+v_r^2}^5}
	+ 3\chi^2 \frac{u(\mu-u)v_r^2}{(1+v_r^2)^3} +\chi \frac{u_r^2}{\sqrt{u^2+u_r^2} \cdot \sqrt{1+v_r^2}^3}
	-\chi^2 \frac{u_r v_r}{(1+v_r^2)^2} \\[1mm]
	\hspace*{19mm}
	+3\chi \cdot \frac{n-1}{r} \cdot \frac{uu_r v_r^2}{\sqrt{u^2+u_r^2} \cdot \sqrt{1+v_r^2}^5}
	- 3\chi^2 \cdot \frac{n-1}{r} \cdot \frac{uv_r^3}{(1+v_r^2)^3}
	\end{array} \right.
  \ee
  for $r\in (0,R)$ and $t\in (0,\tm)$.
\end{lem}
\proof
  We divide (\ref{36.1}) by $u$ and differentiate each term on the right-hand side of the resulting identity separately.
  Using that $u_t=uz$ and hence $u_{rt}=uz_r+u_rz$ and $u_{rrt}=uz_{rr}+2u_r z_r + u_{rr}z$, we first obtain
  \bas
	\bigg(\frac{u^2 u_{rr}}{\sqrt{u^2+u_r^2}^3}\bigg)_t
	&=& \frac{u^2 u_{rrt}}{\sqrt{u^2+u_r^2}^3}
	+ 2\frac{uu_t u_{rr}}{\sqrt{u^2+u_r^2}^3}
	- \frac{3}{2} \cdot \frac{u^2 u_{rr} \cdot (2uu_t + 2u_r u_{rt})}{\sqrt{u^2+u_r^2}^5} \\
	&=& \frac{u^3}{\sqrt{u^2+u_r^2}^3} \cdot z_{rr}
	+ 2\frac{u^2 u_r}{\sqrt{u^2+u_r^2}^3} \cdot z_r
	+ \frac{u^2 u_{rr}}{\sqrt{u^2+u_r^2}^3} \cdot z
	+ 2\frac{u^2 u_{rr}}{\sqrt{u^2+u_r^2}^3} \cdot z \\
	& & - 3\frac{u^4 u_{rr}}{\sqrt{u^2+u_r^2}^5} \cdot z
	- 3\frac{u^3 u_r u_{rr}}{\sqrt{u^2+u_r^2}^5} \cdot z_r
	- 3\frac{u^2 u_r^2 u_{rr}}{\sqrt{u^2+u_r^2}^5} \cdot z.
  \eas
  Since
  \bas
	& & \hspace*{-30mm}
	\frac{u^2 u_{rr}}{\sqrt{u^2+u_r^2}^3} \cdot z
	+ 2\frac{u^2 u_{rr}}{\sqrt{u^2+u_r^2}^3} \cdot z
	- 3\frac{u^4 u_{rr}}{\sqrt{u^2+u_r^2}^5} \cdot z
	- 3\frac{u^2 u_r^2 u_{rr}}{\sqrt{u^2+u_r^2}^5} \cdot z \\
	&=& \frac{u^2 u_{rr}}{\sqrt{u^2+u_r^2}^5} \cdot \Big\{ (u^2+u_r^2)+2(u^2+u_r^2)-3u^2-3u_r^2\Big\} \\[2mm]
	&=& 0,
  \eas
  this yields
  \be{38.3}
	\bigg(\frac{u^2 u_{rr}}{\sqrt{u^2+u_r^2}^3}\bigg)_t
	= \frac{u^3}{\sqrt{u^2+u_r^2}^3} \cdot z_{rr}
	+ \bigg\{ 2\frac{u^2 u_r}{\sqrt{u^2+u_r^2}^3} - 3\frac{u^3 u_r u_{rr}}{\sqrt{u^2+u_r^2}^5} \bigg\} \cdot z_r
  \ee
  for all $r\in (0,R)$ and $t\in (0,\tm)$.\\
  Next,
  \bas
	\bigg( \frac{u_r^4}{u\sqrt{u^2+u_r^2}^3}\bigg)_t
	&=& 4\frac{u_r^3 u_{rt}}{u\sqrt{u^2+u_r^2}^3}
	- \frac{u_r^4 u_t}{u^2 \sqrt{u^2+u_r^2}^3}
	-\frac{3}{2} \cdot \frac{u_r^4 \cdot (2uu_t+2u_r u_{rt})}{u\sqrt{u^2+u_r^2}^5} \\
	&=& 4\frac{u_r^3}{\sqrt{u^2+u_r^2}^3} \cdot z_r
	+ 4\frac{u_r^4}{u\sqrt{u^2+u_r^2}^3} \cdot z
	- \frac{u_r^4}{u\sqrt{u^2+u_r^2}^3} \cdot z \\
	& & -3\frac{uu_r^4}{\sqrt{u^2+u_r^2}^5} \cdot z
	- 3\frac{u_r^5}{\sqrt{u^2+u_r^2}^5} \cdot z_r
	- 3\frac{u_r^6}{u\sqrt{u^2+u_r^2}^5} \cdot z,
  \eas
  where again the zero-order terms have a vanishing sum in the sense that
  \bas
	& & \hspace*{-30mm}
	4\frac{u_r^4}{u\sqrt{u^2+u_r^2}^3} \cdot z
	- \frac{u_r^4}{u\sqrt{u^2+u_r^2}^3} \cdot z
	-3\frac{uu_r^4}{\sqrt{u^2+u_r^2}^5} \cdot z
	- 3\frac{u_r^6}{u\sqrt{u^2+u_r^2}^5} \cdot z \\
	&=& \frac{u_r^4}{u\sqrt{u^2+u_r^2}^5} \cdot \Big\{ 4(u^2+u_r^2) - (u^2+u_r^2) - 3u^2 - 3u_r^2 \Big\} \cdot z \\[2mm]
	&=& 0,
  \eas
  so that
  \be{38.4}
	\bigg(\frac{u_r^4}{u\sqrt{u^2+u_r^2}^3}\bigg)_t
	= \bigg\{ 4\frac{u_r^3}{\sqrt{u^2+u_r^2}^3} - 3\frac{u_r^5}{\sqrt{u^2+u_r^2}^5} \bigg\} \cdot z_r
  \ee
  for all $r\in (0,R)$ and $t\in (0,\tm)$.\\
  Likewise,
  \bas
	\bigg(\frac{u_r}{\sqrt{u^2+u_r^2}}\bigg)_t
	&=& \frac{u_{rt}}{\sqrt{u^2+u_r^2}} - \frac{1}{2} \cdot \frac{u_r\cdot (2uu_t+2u_r u_{rt})}{\sqrt{u^2+u_r^2}^3} \\
	&=& \frac{u}{\sqrt{u^2+u_r^2}} \cdot z_r + \frac{u_r}{\sqrt{u^2+u_r^2}} \cdot z
	- \frac{u^2 u_r}{\sqrt{u^2+u_r^2}^3} \cdot z
	- \frac{uu_r^2}{\sqrt{u^2+u_r^2}^3} \cdot z_r
	- \frac{u_r^3}{\sqrt{u^2+u_r^2}^3} \cdot z \\
	&=& \frac{u}{\sqrt{u^2+u_r^2}^3} \cdot \Big\{ (u^2+u_r^2)-u_r^2\Big\} \cdot z_r \\
	& & + \frac{u_r}{\sqrt{u^2+u_r^2}^3} \cdot \Big\{(u^2+u_r^2)-u^2-u_r^2\Big\} \cdot z \\
	&=& \frac{u^3}{\sqrt{u^2+u_r^2}^3} \cdot z_r,
  \eas
  whence
  \be{38.5}
	\bigg(\frac{n-1}{r} \frac{u_r}{\sqrt{u^2+u_r^2}} \bigg)_t
	= \frac{n-1}{r} \cdot \frac{u^3}{\sqrt{u^2+u_r^2}^3} \cdot z_r
  \ee
  for $r\in (0,R)$ and $t\in (0,\tm)$.\\
  As for the respective terms originating from the rightmost three summands in (\ref{36.1}), we make use of (\ref{33.3})
  to express $v_{rt}$ conveniently. We thereby compute
  \bea{38.6}
	\bigg(-\chi \frac{\mu-u}{\sqrt{1+v_r^2}^3}\bigg)_t
	&=& \chi\frac{u_t}{\sqrt{1+v_r^2}^3} + \frac{3}{2} \chi \frac{\mu-u}{\sqrt{1+v_r^2}^5} \cdot 2v_r v_{rt} \nn\\
	&=& \chi\frac{u}{\sqrt{1+v_r^2}^3} \cdot z
	+ 3\chi \frac{(\mu-u)v_r}{\sqrt{1+v_r^2}^5} \cdot
	\Big\{ -\frac{uu_r}{\sqrt{u^2+u_r^2}} + \chi \frac{uv_r}{\sqrt{1+v_r^2}} \Big\} \nn\\
	&=& \chi\frac{u}{\sqrt{1+v_r^2}^3} \cdot z
	- 3\chi \frac{u(\mu-u)u_r v_r}{\sqrt{u^2+u_r^2} \cdot \sqrt{1+v_r^2}^5}
	+ 3\chi^2 \frac{u(\mu-u)v_r^2}{(1+v_r^2)^3}
  \eea
  and
  \bea{38.7}
	\bigg( -\chi \frac{u_r v_r}{u\sqrt{1+v_r^2}}\bigg)_t
	&=& - \chi \frac{u_{rt} v_r}{u\sqrt{1+v_r^2}} - \chi \frac{u_r v_{rt}}{u\sqrt{1+v_r^2}}
	+ \chi \frac{u_r v_r u_t}{u^2 \sqrt{1+v_r^2}}
	+ \frac{1}{2} \chi \frac{u_r v_r}{u\sqrt{1+v_r^2}^3} \cdot 2v_r v_{rt} \nn\\
	&=& - \chi \frac{v_r}{\sqrt{1+v_r^2}} \cdot z_r
	- \chi \frac{u_r v_r}{u\sqrt{1+v_r^2}} \cdot z \nn\\
	& & - \chi \frac{u_r}{u\sqrt{1+v_r^2}} \cdot
	\Big\{ -\frac{uu_r}{\sqrt{u^2+u_r^2}} + \chi \frac{uv_r}{\sqrt{1+v_r^2}} \Big\} \nn\\
	& & + \chi \frac{u_r v_r}{u\sqrt{1+v_r^2}} \cdot z
	+ \chi \frac{u_r v_r^2}{u\sqrt{1+v_r^2}^3} \cdot
	\Big\{ -\frac{uu_r}{\sqrt{u^2+u_r^2}} + \chi \frac{uv_r}{\sqrt{1+v_r^2}} \Big\} \nn\\
	&=& - \chi \frac{v_r}{\sqrt{1+v_r^2}} \cdot z_r
	+ \chi \frac{u_r^2}{\sqrt{u^2+u_r^2} \cdot \sqrt{1+v_r^2}}
	-\chi^2 \frac{u_r v_r}{(1+v_r^2)^2} \nn\\
	&=& - \chi \frac{v_r}{\sqrt{1+v_r^2}} \cdot z_r
	+ \chi \frac{u_r^2}{\sqrt{u^2+u_r^2} \cdot \sqrt{1+v_r^2}^3}
	-\chi^2 \frac{u_r v_r}{(1+v_r^2)^2},
  \eea
  and observe that
  \bas
	\bigg(\frac{v_r^3}{\sqrt{1+v_r^2}^3}\bigg)_t
	&=& \Big\{ \frac{3v_r^2}{\sqrt{1+v_r^2}^3} - \frac{3}{2}
	\cdot \frac{v_r^3 \cdot 2v_r}{\sqrt{1+v_r^2}^5} \Big\} \cdot v_{rt} \nn\\
	&=& - 3\frac{uu_r v_r^2}{\sqrt{u^2+u_r^2} \cdot \sqrt{1+v_r^2}^5}
	+3\chi \frac{uv_r^3}{(1+v_r^2)^3}
  \eas
  to obtain
  \be{38.8}
	\bigg( -\chi \cdot \frac{n-1}{r} \cdot \frac{v_r^3}{\sqrt{1+v_r^2}^3} \bigg)_t
	= 3\chi\cdot \frac{n-1}{r} \cdot \frac{uu_r v_r^2}{\sqrt{u^2+u_r^2} \cdot \sqrt{1+v_r^2}^5}
	-3\chi^2 \cdot \frac{n-1}{r} \cdot \frac{uv_r^3}{(1+v_r^2)^3}
  \ee
  for $r\in (0,R)$ and $t\in (0,\tm)$.
  In light of (\ref{38.3})-(\ref{38.8}), (\ref{36.1}) easily yields (\ref{38.1}) with $B_1, B_{21}, B_{22}, B_3$
  and $B_4$ as in (\ref{38.2}).
\qed
On suitably estimating the inhomogeneous term $B_4$ herein by means of Corollary \ref{cor40},
we can develop (\ref{38.1}) into a nonlocal parabolic inequality for $z$ as follows.
\begin{lem}\label{lem39}
  Suppose that $\tm<\infty$, but that $\sup_{(r,t)\in (0,R)\times (0,\tm)} u(r,t)<\infty$.
  Then there exist a constant $d>0$ and continuous functions $b_1, b_{21}, b_{22}$ and $b_3$ on $[0,R]\times [0,\tm)$
  with the properties that $b_1$ and $b_{22}$ are nonnegative, that $b_3$ is bounded on $(0,R)\times (0,\tm)$, and such
  that $z=\frac{u_t}{u}$ satisfies
  \be{39.1}
	z_t(r,t) \le b_1(r,t) z_{rr} + b_{21}(r,t) z_r + \frac{b_{22}(r,t)}{r} z_r + b_3(r,t)z
	+ d \cdot \Big( 1+\|z_+\|_{L^\infty((0,R)\times (0,t))} \Big)
  \ee
  for all $r\in (0,R)$ and $t\in (0,\tm)$.
\end{lem}
\proof
  With $B_1, B_{21}, B_{22}, B_3$ and $B_4$ taken from Lemma \ref{lem38}, we let $b_1:=B_1, b_{21}:=B_{21}, b_{22}:=B_{22}$
  and $b_3:=B_3$. Then from (\ref{38.2}) we immediately obtain that $b_1, b_{21}, b_{22}$ and $b_3$ are continuous in
  $[0,R] \times [0,\tm)$, and that $b_1 \ge 0$ and $b_{22} \ge 0$. Since our boundedness assumption on $u$ ensures that $b_3$
  is bounded, it remains to control the inhomogeneity $B_4$ in (\ref{38.1}) adequately. To this end, we once more
  use our hypothesis along with Lemma \ref{lem333} to pick positive constants $c_1, c_2$ and $c_3$ such that
  \be{39.2}
	u(r,t) \le c_1, \quad
	|v_r(r,t)| \le c_2 r
	\quad \mbox{and} \quad
	|v_{rr}(r,t)| \le c_3
	\qquad \mbox{for all $r\in (0,R)$ and } t\in (0,\tm).
  \ee
  Then in (\ref{38.2}) we can estimate
  \bas
	-3\chi \frac{u(\mu-u) u_r v_r}{\sqrt{u^2+u_r^2} \cdot \sqrt{1+v_r^2}^5}
	&\le& 3\chi \cdot c_1 (\mu+c_1) \cdot c_2 R \cdot \frac{|u_r|}{\sqrt{u^2+u_r^2}}  \\
	&\le& 3c_1(\mu+c_1)c_2 \chi R
  \eas
  and
  \bas
	3\chi^2 \frac{u(\mu-u) v_r v_{rr}}{(1+v_r^2)^3}
	&\le& 3\chi^2 c_1 (\mu+c_1) \cdot c_2 R \cdot c_3 \\
	&\le& 3c_1(\mu+c_1)c_2 c_3 \chi^2 R
  \eas
  as well as
  \bas
	3\chi \cdot \frac{n-1}{r} \cdot \frac{uu_r v_r^2}{\sqrt{u^2+u_r^2} \cdot \sqrt{1+v_r^2}^5}
	&\le& 3\chi \cdot \frac{n-1}{r} \cdot c_1 \cdot c_2^2 r^2 \cdot \frac{|u_r|}{\sqrt{u^2+u_r^2}} \\
	&\le& 3(n-1) c_1 c_2^2 \chi R
  \eas
  and
  \bas
	-3\chi^2 \cdot \frac{n-1}{r} \cdot \frac{uv_r^3}{(1+v_r^2)^3}
	&\le& 3\chi^2 \cdot \frac{n-1}{r} \cdot c_1 \cdot c_2^3 r^3 \\
	&\le& 3(n-1) c_1 c_2^3 \chi^2 R^2
  \eas
  for all $r\in (0,R)$ and $t\in (0,\tm)$.
  In the third and fourth summands in the definition (\ref{38.2}) of $B_4$, however, apparently we can only estimate
  \bas
	\chi \frac{u_r^2}{\sqrt{u^2+u_r^2} \cdot \sqrt{1+v_r^2}^3}
	\le \chi \frac{u_r^2}{\sqrt{u^2+u_r^2}} \le \chi |u_r|
  \eas
  and
  \bas
	- \chi^2 \frac{u_r v_r}{(1+v_r^2)^2} \le \chi^2 |u_r|
  \eas
  for all $r\in (0,R)$ and $t\in (0,\tm)$, with the possibly unbounded factors $|u_r|$ remaining.
  Fortunately, applying Corollary \ref{cor40} yields $c_4>0$ such that
  \bas
	|u_r(r,t)| \le c_4 \cdot \Big(1+\|z_+\|_{L^\infty((0,R)\times (0,t))} \Big)
	\qquad \mbox{for all $r\in (0,R)$ and } t\in (0,\tm).
  \eas
  Therefore, (\ref{39.1}) results from (\ref{38.1}) if we choose $d>0$ conveniently large.
\qed
\subsection{Boundedness of $z$ from above}
Apparently, nonlocal parabolic inequalities of type (\ref{39.1}) do not allow for general comparison priciples.
After all, the fact that here the memory term enjoys a certain linear boundedness property with respect to $z_+$
enables us to follow a maximum principle-type reasoning to establish an essentially exponential upper bound for $z$
and thereby obtain the following.
\begin{lem}\label{lem41}
  Assume that $\tm<\infty$ and $\sup_{(r,t)\in (0,R)\times (0,\tm)} u(r,t)<\infty$.
  Then there exists $C>0$ such that $z=\frac{u_t}{u}$ satisfies
  \be{41.1}
	z(r,t) \le C
	\qquad \mbox{for all $r\in (0,R)$ and } t\in (0,\tm).
  \ee
\end{lem}
\proof
  We let $b_1, b_{21}, b_{22}, b_3$ and $d$ be as provided by Lemma \ref{lem39}, so that by boundedness of $b_3$ we can
  find $c_1>0$ such that
  \be{41.2}
	b_3(r,t) \le c_1
	\qquad \mbox{for all $r\in (0,R)$ and } t\in (0,\tm).
  \ee
  We than fix $\alpha>0$ large enough fulfilling
  \be{41.3}
	\alpha>c_1+d
  \ee
  and let
  \bas
	\varphi(r,t):=e^{-\alpha t} z(r,t) - dt
	\qquad \mbox{for $r\in [0,R]$ and } t\in [0,\tm).
  \eas
  Then according to Lemma \ref{lem39},
  \bea{41.4}
	\varphi_t
	&=& e^{-\alpha t} (z_t-\alpha z) -d \nn\\
	&\le& e^{-\alpha t} \cdot \Big\{ b_1(r,t) z_{rr} + b_{21}(r,t) z_r + \frac{b_{22}}{r} z_r
	+ b_3(r,t)z + d\|z_+\|_{L^\infty((0,R)\times (0,t))} + d - \alpha z \Big\} \, -d \nn\\
	&=& b_1(r,t) \varphi_{rr} + b_{21}(r,t)\varphi_r + \frac{b_{22}(r,t)}{r} \varphi_r
	+ \Big(b_3(r,t)-\alpha\Big) \cdot \Big(\varphi+dt \Big) \nn\\[1mm]
	& & + d e^{-\alpha t} \|z_+\|_{L^\infty((0,R)\times (0,t))}
	+ d e^{-\alpha t} - d
	\qquad \mbox{for all $r\in (0,R)$ and } t\in (0,\tm),
  \eea
  and since
  $z_r=\big(\frac{u_t}{u}\big)_r = \frac{u_{rt}}{u}-\frac{u_r u_t}{u^2}$ in $[0,R] \times [0,\tm)$,
  the fact that $u_r(0,t)=u_r(R,t)=0$ for all $t\in (0,\tm)$ entails that	
  \be{41.5}
	\varphi_r(0,t)=\varphi_r(R,t)=0
	\qquad \mbox{for all } t\in (0,\tm).
  \ee
  Now if for some $T\in (0,\tm)$, the value $S:=\sup_{(r,t)\in (0,R)\times (0,T)} \varphi(r,t)$ was positive
  and attained at some point $(r_0,t_0)\in [0,R]\times [0,T]$ with $t_0>0$, then necessarily
  \be{41.6}
	\varphi_t(r_0,t_0) \ge 0,
  \ee
  and (\ref{41.5}) ensures that in both cases $r_0\in (0,R)$ and $r_0\in \{0,R\}$ we moreover must have
  \be{41.7}
	\varphi_r(r_0,t_0)=0
	\qquad \mbox{and} \varphi_{rr}(r_0,t_0) \le 0.
  \ee
  We claim that these properties imply that
  \be{41.8}
	0 \le \Big(b_3(r_0,t_0)-\alpha\Big) \cdot \Big( \varphi(r_0,t_0)+dt_0\Big)
	+ d e^{-\alpha t_0} \|z_+\|_{L^\infty((0,R)\times (0,t_0))} + d e^{-\alpha t_0} - d.
  \ee
  Indeed, in the case $r_0\in (0,R)$ we may directly apply (\ref{41.4}) to easily deduce this from (\ref{41.6})
  and (\ref{41.7}).
  When $r_0=R$, by continuity of $\varphi, \varphi_t, \varphi_r$ and $\varphi_{rr}$ in $[0,R]\times (0,\tm)$
  it is clear that (\ref{41.4}) actually remains valid at $(r_0,t_0)$, so that (\ref{41.8}) follows from the same argument.
  If $r_0=0$, however, we make use of the favorable sign of the singular term $\frac{b_{22}}{r}$ in (\ref{41.4})
  by first
  choosing, once more relying on the extremal property of $\varphi(r_0,t_0)$, a sequence $(r_j)_{j\in\N} \subset (0,R)$
  such that $r_j\searrow 0$ as $j\to\infty$ and
  \bas
	\varphi_r(r_j,t_0) \le 0
	\qquad \mbox{for all } j\in\N,
  \eas
  and then evaluating (\ref{41.4}) at $r=r_j$ to see that
  \bas
	\varphi_t(r_j,t_0)
	&\le& b_1(r_j,t_0) \varphi_{rr}(r_j,t_0) + b_{21}(r_j,t_0) \varphi_r(r_j,t_0)
	+ \Big(b_3(r_j,t_0)-\alpha \Big) \cdot \Big( \varphi (r_j,t_0) + dt_0 \Big) \\
	& & + d e^{-\alpha t_0} \|z_+\|_{((0,R) \times (0,t_0))}
	+ d e^{-\alpha t_0} -d
  \eas
  for all $j\in\N$. Again by continuity of $\varphi,\varphi_t, \varphi_r$ and $\varphi_{rr}$, we may take $j\to\infty$
  to conclude that
  \bas
	\varphi_t(0,t_0)
	&\le& b_1(0,t_0) \varphi_{rr}(0,t_0) + b_{21}(0,t_0) \varphi_r(0,t_0)
	+ \Big( b_3(0,t_0) -\alpha \Big) \cdot \Big( \varphi(0,t_0) + dt_0\Big) \\
	& & + d e^{-\alpha t_0} \|z_+\|_{((0,R)\times (0,t_0))}
	+ d e^{-\alpha t_0} - d,
  \eas
  whereupon one more application of (\ref{41.6}) and (\ref{41.7}) yields (\ref{41.8}) also in this case.\abs
  Now observing that $e^{-\alpha t_0} \le 1$ and using that $S=\varphi(r_0,t_0)$ is positive, in view of (\ref{41.2})
  we obtain from (\ref{41.8}) that
  \be{41.9}
	0 \le (c_1-\alpha) \cdot \Big( \varphi(r_0,t_0)+dt_0\Big)
	+ d e^{-\alpha t_0} \|z_+\|_{L^\infty((0,R)\times (0,t_0))}.
  \ee
  Here we rewrite $z= e^{\alpha t} \varphi + dt e^{\alpha t}$ and use that if $f$ and $g$ are two functions on a set
  $D\subset \R^N$, $N\ge 1$, then both $(f+g)_+$ and $\sup_D \{f+g\} \le \sup_D f + \sup_D g$. We thereby obtain that
  \bas
	d e^{-\alpha t_0} \|z_+\|_{L^\infty((0,R)\times (0,t_0))}
	&=& d e^{-\alpha t_0} \cdot
	\sup_{(r,s)\in (0,R)\times (0,t_0)} \Big\{ e^{\alpha s}\varphi(r,s) + ds e^{\alpha s} \Big\} \\
	&\le& d e^{\alpha t_0} \cdot
	\sup_{(r,s)\in (0,R)\times (0,t_0)} \Big\{ e^{\alpha s}\varphi_+(r,s) + ds e^{\alpha s} \Big\} \\
	&\le& d e^{\alpha t_0} \cdot \bigg\{
	\sup_{(r,s)\in (0,R)\times (0,t_0)} \Big\{ e^{\alpha s}\varphi_+(r,s) \Big\}
	+ \sup_{s \in (0,t_0)} \Big\{ ds e^{\alpha s} \Big\} \bigg\}\\
	&\le& d e^{\alpha t_0} \cdot
	\sup_{(r,s)\in (0,R)\times (0,t_0)} \Big\{ e^{\alpha s}\varphi_+(r,s) \Big\}
	+ d^2 t_0.
  \eas
  Since from the definition of $S$ we know that
  \bas
	\sup_{(r,s)\in (0,R)\times (0,t_0)} \Big\{ e^{\alpha s}\varphi_+(r,s) \Big\}
	\le e^{\alpha t_0} \cdot
	\sup_{(r,s)\in (0,R)\times (0,t_0)} \varphi_+(r,s)
	= e^{\alpha t_0} \cdot \varphi(r_0,t_0),
  \eas
  this entails that
  \bas
	d e^{-\alpha t_0} \|z_+\|_{L^\infty((0,R)\times (0,t_0))}
	\le d e^{-\alpha t_0} \|z_+\|_{L^\infty((0,R)\times (0,t_0))}
	\le d \varphi(r_0,t_0) + d^2 t_0,
  \eas
  so that (\ref{41.9}) yields the inequality
  \bas
	0 &\le& (c_1-\alpha)\cdot \Big( \varphi(r_0,t_0)+dt_0 \Big) + d\varphi(r_0,t_0) + dt_0 \\
	&=& (c_1-\alpha+d) \cdot \Big(\varphi(r_0,t_0)+dt_0\Big).
  \eas
  In light of our restriction (\ref{41.3}) on $\alpha$, however, this contradicts the positivity of $\varphi(r_0,t_0)$
  and thereby proves that actually $\varphi$ cannot attain a positive maximum over any such region $[0,R]\times [0,T]$,
  $T\in (0,\tm)$, at a positive time $t_0$. This means that in fact
  \bas
	\varphi(r,t) \le \|\varphi_+(\cdot,0)\|_{L^\infty((0,R))} = \|z_+(\cdot,0)\|_{L^\infty((0,R))}
	\qquad \mbox{for all $r\in (0,R)$ and } t\in (0,\tm)
  \eas
  and hence
  \bas
	z(r,t) &=& e^{\alpha t} \Big(\varphi(r,t) + dt \Big) \\
	&\le& e^{\alpha \tm} \cdot \Big\{ \|z_+(\cdot,0)\|_{L^\infty((0,R))} + d\tm\Big\}
	\qquad \mbox{for all $r\in (0,R)$ and } t\in (0,\tm),
  \eas
  which establishes (\ref{41.1}).
\qed
\subsection{Boundedness of $u$ implies extensibility. Proof of Theorem \ref{theo43}}
Combining the latter lemma with Corollary \ref{cor40} now directly yields the desired bound for $u_r$.
\begin{cor}\label{cor42}
  If $\tm<\infty$ but $\sup_{(r,t)\in (0,R)\times (0,\tm)} u(r,t)<\infty$,
  then there exists $C>0$ such that
  \bas
	\|u_r(\cdot,t)\|_{L^\infty((0,R))} \le C
	\qquad \mbox{for all } t\in (0,\tm).
  \eas
\end{cor}
\proof
  Thanks to the upper estimate for $z$ obtained in Lemma \ref{lem41}, this is an immediate consequence of Corollary
  \ref{cor40}.
\qed
We can thereby readily verify our main statement on local existence and extensibility.\abs
\proofc of Theorem \ref{theo43}.\quad
  In view of the local existence result established in Lemma \ref{lem21}, we only need to verify (\ref{43.1}).
  Indeed, if (\ref{43.1}) was false, then for some solution the respective maximal
  existence time would satisfy $\tm<\infty$ but $\limsup_{t\nearrow\tm} \|u(\cdot,t)\|_{L^\infty(\Omega)}<\infty$.
  Then, however, Corollary \ref{cor42} would apply to assert that also
  $\limsup_{t\nearrow\tm} \|\nabla u(\cdot,t)\|_{L^\infty(\Omega)}$ would be finite.
  Along with the lower bound for $u$ provided by Lemma \ref{lem23}, this would contradict the extensibility criterion
  (\ref{ext_crit}) in Lemma \ref{lem21}.
\qed
\mysection{Boundedness for small $\chi$. Proof of Theorem \ref{theo55}}\label{sect6}
In light of the extensibility criterion provided by Theorem \ref{theo43}, in order to prove both global existence and
boundedness of a solution it is sufficient to derive an a priori estimate for $(u,v)$ in $(L^\infty(\Omega\times (0,T)))^2$
which does not explicitly depend on $T<\tm \le \infty$.
As a preparation for the proof of this in Theorem \ref{theo55} below, let us state the following elementary inequality.
\begin{lem}\label{lem51}
  Let $p\ge 1$. Then
  \be{51.1}
	\io u^{p-1}|\nabla u| \le \io \frac{u^{p-1}|\nabla u|^2}{\sqrt{u^2+|\nabla u|^2}}
	+ \io u^p
	\qquad \mbox{for all } t\in (0,\tm).
  \ee
\end{lem}
\proof
  By means of Young's inequality, we see that
  \be{51.2}
	\io u^{p-1}|\nabla u|
	\le \frac{1}{2} \io \frac{u^{p-1}|\nabla u|^2}{\sqrt{u^2+|\nabla u|^2}}
	+ \frac{1}{2} \io u^{p-1}\sqrt{u^2+|\nabla u|^2}
	\qquad \mbox{for all } t\in (0,\tm),
  \ee
  where using the elementary inequality $\sqrt{X+Y} \le \sqrt{X}+\sqrt{Y}$, valid for all $X\ge 0$ and $Y\ge 0$,
  we can estimate
  \bas
	\frac{1}{2} \io u^{p-1}\sqrt{u^2+|\nabla u|^2}
	\le \frac{1}{2} \io u^p + \frac{1}{2} \io u^{p-1}|\nabla u|.
  \eas
  Therefore, (\ref{51.1}) results from (\ref{51.2}).
\qed
We are now in the position to make sure that if either $n\ge 2, \chi<1$ and $u_0$ is an arbitrary function satisfying
(\ref{init}), or $n=1$, $\chi>0$ and $\io u_0<m_c$ with $m_c$ as in (\ref{mc}),
then the solution of (\ref{0}) in fact is global and remains bounded:\abs
\proofc of Theorem \ref{theo55}.\quad
  We let $p_k:=2^k$ and, given $T\in (0,\tm)$, introduce
  \be{55.22}
	M_k:=\sup_{t\in (0,T)} \io u^{p_k}(x,t)dx
  \ee
  for nonnegative integers $k$.
  Then clearly $M_k$ is well-defined for any such $T$ and $k$, and in order to control $M_k$ appropriately,
  we fix $k\ge 1$ and multiply the first equation in (\ref{0}) by $pu^{p-1}$ for $p:=p_k$ to see upon integrating
  by parts that
  \bea{55.2}
	\frac{d}{dt} \io u^p
	+ p(p-1) \io \frac{u^{p-1}|\nabla u|^2}{\sqrt{u^2+|\nabla u|^2}}
	&=& p(p-1)\chi \io \frac{u^{p-1}\nabla u\cdot \nabla v}{\sqrt{1+|\nabla v|^2}} \nn\\
	&\le& p(p-1)\chi \io u^{p-1}|\nabla u| \cdot \frac{|\nabla v|}{\sqrt{1+|\nabla v|^2}}
  \eea
  for all $t\in (0,\tm)$.
  Here in the multi-dimensional case, in which no evident uniform a priori bound for $|\nabla v|$ seems available,
  we use the trivial pointwise inequality $\frac{|\nabla v|}{\sqrt{1+|\nabla v|^2}} \le 1$ to obtain
  \be{55.211}
	\frac{d}{dt} \io u^p
	+ p(p-1) \io \frac{u^{p-1}|\nabla u|^2}{\sqrt{u^2+|\nabla u|^2}}
	\le p(p-1)\chi \io u^{p-1}|\nabla u|
	\quad \mbox{for all } t\in (0,\tm)
	\qquad \mbox{if } n\ge 2.
  \ee
  In the one-dimensional setting, however, from (\ref{333.0}) and (\ref{mu}) we know that $|\nabla v|=|v_r| \le m$ throughout
  $\Omega\times (0,\tm)$, whence by monotonicity of $0 \le \xi \mapsto \frac{\xi}{\sqrt{1+\xi^2}}$ we infer from
  (\ref{55.2}) that
  \be{55.212}
	\frac{d}{dt} \io u^p
	+ p(p-1) \io \frac{u^{p-1}|\nabla u|^2}{\sqrt{u^2+|\nabla u|^2}}
	\le p(p-1)\chi \cdot \frac{m}{\sqrt{1+m^2}} \cdot \io u^{p-1}|\nabla u|
	\quad \mbox{for all } t\in (0,\tm)
	\quad \mbox{if } n=1.
  \ee
  In both (\ref{55.211}) and (\ref{55.212}) we now apply Lemma \ref{lem51} to estimate
  \bas
	p(p-1) \io \frac{u^{p-1}|\nabla u|^2}{\sqrt{u^2+|\nabla u|^2}}
	\ge p(p-1) \io u^{p-1}|\nabla u| - p(p-1) \io u^p
  \eas
  and thus obtain on writing
  \be{Lambda}
	\Lambda:=\left\{ \begin{array}{ll}
	\frac{m}{\sqrt{1+m^2}} \qquad & \mbox{if } n=1, \\[1mm]
	1 & \mbox{if } n\ge 2,
	\end{array} \right.
  \ee
  and adding $\io u^p$ on both sides of (\ref{55.211}) and (\ref{55.212}), respectively, that
  \bea{55.3}
	\frac{d}{dt} \io u^p + \io u^p + p(p-1)(1-\chi \Lambda) \io u^{p-1}|\nabla u|
	&\le& \Big\{p(p-1)+1\Big\} \cdot \io u^p \nn\\
	&\le& p^2 \io u^p
	\qquad \mbox{for all } t\in (0,\tm).
  \eea
  We next invoke the Gagliardo-Nirenberg inequality \cite{win_critexp}
  to find $c_1>0$ such that with $a:=\frac{n}{n+1}$ we have
  \bas
	\|\varphi\|_{L^1(\Omega)}
	\le c_1 \|\nabla\varphi\|_{L^1(\Omega)}^a \|\varphi\|_{L^\frac{1}{2}(\Omega)}^{1-a}
	+ c_1\|\varphi\|_{L^\frac{1}{2}(\Omega)}
	\qquad \mbox{for all } \varphi\in W^{1,1}(\Omega),
  \eas
  and thereby obtain that
  \bas
	\io u^p \le c_1 \Big( \io |\nabla u^p| \Big)^a \cdot \Big(\io u^\frac{p}{2}\Big)^{2(1-a)}
	\qquad \mbox{for all } t\in (0,\tm).
  \eas
  Since our specification of $p=p_k=2^k$ allows us to use the definition (\ref{55.22}) of $M_{k-1}$ in estimating
  \bas
	\io u^\frac{p}{2} \le M_{k-1}
	\qquad \mbox{for all } t\in (0,T),
  \eas
  this implies that
  \bas
	\io u^p \le c_1 \Big( \io |\nabla u^p|\Big)^a \cdot M_{k-1}^{2(1-a)} + c_1 M_{k-1}^2
	\qquad \mbox{for all } t\in (0,T).
  \eas
  Thanks to the fact that our assumptions ensure that $\chi\Lambda<1$, another application of Young's inequality therefore
  provides $c_2>0$ fulfilling
  \bas
	p^2 \io u^p
	\le (p-1)(1-\chi\Lambda) \io |\nabla u^p|
	+ c_2 p^\frac{2}{1-a} M_{k-1}^2 + c_1 p^2 M_{k-1}^2
	\qquad \mbox{for all } t\in (0,T),
  \eas
  from which due to the evident fact that $p^2 \le p^\frac{2}{1-a}$ we obtain that
  \bas
	p^2 \io u^p
	\le p(p-1)(1-\chi\Lambda) \io u^{p-1}|\nabla u|
	+ c_3 p^\frac{2}{1-a} M_{k-1}^2
	\qquad \mbox{for all } t\in (0,T)
  \eas
  with $c_3:=c_1+c_2$.
  Therefore, (\ref{55.3}) entails the autonomous ODI
  \bas
	\frac{d}{dt} \io u^p + \io u^p
	\le c_3 p^\frac{2}{1-a} M_{k-1}^2,
	\qquad t\in (0,T),
  \eas
  for $(0,T) \ni t \mapsto \io u^p(x,t)dx$, which upon a comparison argument implies that
  \be{55.4}
	M_k \le \max \bigg\{ \io u_0^{p_k} \, , \, c_3 p_k^\frac{2}{1-a} M_{k-1}^2 \bigg\}
	\qquad \mbox{for all } k\ge 1.
  \ee
  Now if there exists a sequence $(k_j)_{j\in\N} \subset \N$ such that $k_j\to\infty$ as $j\to\infty$ and
  \be{55.44}
	M_{k_j} \le \io u_0^{p_{k_j}}
	\qquad \mbox{for all } j\in\N,
  \ee
  we may take the $p_{k_j}$-th root on both sides here to see that according to the definition (\ref{55.22})
  of $M_{k_j}$ we have
  \bas
	\sup_{t\in (0,T)} \|u(\cdot,t)\|_{L^{p_{k_j}}(\Omega)} \le \|u_0\|_{L^{p_{k_j}}(\Omega)},
  \eas
  which on letting $j\to\infty$ implies that
  \be{55.5}
	\sup_{t\in (0,T)} \|u(\cdot,t)\|_{L^\infty(\Omega)} \le \|u_0\|_{L^\infty(\Omega)}
  \ee
  in this case.\\
  Conversely, if no such sequence exists, then (\ref{55.4}) means that with some suitably large $k_0\in\N$ we have
  \bas
	M_k \le c_3 p_k^\frac{2}{1-a} M_{k-1}
	\qquad \mbox{for all } k\ge k_0.
  \eas
  Since $p_k^\frac{2}{1-a}=(2^\frac{2}{1-a})^k$, it is easy to see that this entails the existence of a number $b>1$
  independent of $T$ which satisfies
  \be{55.6}
	M_k \le b^k M_{k-1}
	\qquad \mbox{for all } k\ge 1.
  \ee
  By a straightforward induction, this warrants that
  \bas
	M_k \le b^{\sum_{j=0}^k j\cdot 2^{k-j}} \cdot M_0^{2^k}
	\qquad \mbox{for all } k\ge 1,
  \eas
  where by an elementary computation,
  \bas
	\sum_{j=0}^k j\cdot 2^{k-j}
	&=& 2^{k-1} \cdot \sum_{j=0}^k j\cdot \Big(\frac{1}{2}\Big)^{j-1} \\
	&=& 2^{k-1} \cdot \frac{k\cdot(\frac{1}{2})^{k+1} -(k+1)\cdot (\frac{1}{2})^k + 1}{(\frac{1}{2})^2} \\
	&=& k - 2(k+1) + 2^{k+1} \\[2mm]
	&\le& 2^{k+1}
	\qquad \mbox{for all } k\ge 1.
  \eas
  Thus,
  \bas
	M_k \le B^{2^{k+1}} \cdot M_0^{2^k}
	\qquad \mbox{for all } k\ge 1,
  \eas
  by (\ref{55.22}) implying that
  \be{55.7}
	\sup_{t\in (0,T)} \|u(\cdot,t)\|_{L^{p_k}(\Omega)}
	= M_k^\frac{1}{2^k}
	\le b^2 M_0
	\qquad \mbox{for all } k\ge 1.
  \ee
  Now as by the evident mass conservation property in (\ref{0}) we have $\io u(x,t)dx=\io u_0(x)dx$ for all
  $t\in (0,\tm)$ and hence $M_0=\|u_0\|_{L^1(\Omega)}$, taking $k\to\infty$ in (\ref{55.7}) shows that in this
  second case,
  \be{55.8}
	\sup_{t\in (0,T)} \|u(\cdot,t)\|_{L^\infty(\Omega)}
	\le b^2 \|u_0\|_{L^1(\Omega)}.
  \ee
  Since all expressions on the right-hand sides of (\ref{55.5}) and (\ref{55.8}) do not depend on $T\in (0,\tm)$,
  and since boundedness of $u$ clearly implies boundedness of $v$ by standard elliptic estimates, the proof
  is complete.
\qed

\end{document}